\documentclass[12pt]{amsart}
\usepackage{amsmath}
\usepackage{amssymb}
\usepackage{amsthm}
\usepackage{psfig}
\newcommand{\ignore}[1]{\relax}

\newcommand{\C}{\mathbb C}
\newcommand{\R}{\mathbb R}
\newcommand{\Z}{\mathbb Z}
\newcommand{\N}{\mathbb N}

\newcommand{\s}{\mathcal S}

\newcommand{\ind}{\operatorname{ind}}

\newcommand{\re}{\operatorname{Re}}

\newcommand{\Val}{\operatorname{Val}}
\newcommand{\val}{\operatorname{val}}

\newcommand{\Tau}{\mathcal T}

\newcommand{\ncl}{N^{\operatorname{irr}}}
\newcommand{\nclr}{N}

\newtheorem{lem}{Lemma}[section]

\newtheorem{cor}[lem]{Corollary}
\newtheorem{thma}[lem]{Theorem}
\newtheorem{thm}[lem]{Theorem}
\newtheorem{coro}[lem]{Corollary}
\newtheorem{prop}[lem]{Proposition}
\newtheorem{add}[lem]{Addendum}
\theoremstyle{definition}
\newtheorem{defn}[lem]{Definition}
\newtheorem{exa}[lem]{Example}
\newtheorem{que}[lem]{Question}
\theoremstyle{remark}
\newtheorem{rmk}[lem]{Remark}
\newtheorem{rem}[lem]{Remark}

\newcommand{\Rtr}{{\mathbb R}_{\operatorname{trop}}}
\newcommand{\tor}{(\C^*)^{n}}
\newcommand{\tordva}{(\C^*)^{2}}
\newcommand{\tortri}{(\C^*)^{3}}

\newcommand{\rtor}{(\R^*)^{n}}

\newcommand{\mult}{\operatorname{mult}}

\newcommand{\dH}{d_{\operatorname{Haus}}}
\newcommand{\dd}{\partial}
\newcommand{\am}{\mathcal{A}}

\newcommand{\cp}{{\mathbb C}{\mathbb P}}
\newcommand{\rp}{{\mathbb R}{\mathbb P}}

\newcommand{\ppp}{{\mathcal P}}
\newcommand{\qqq}{{\mathcal R}}

\newcommand{\Log}{\operatorname{Log}}

\newcommand{\Vol}{\operatorname{Vol}}

\newcommand{\Int}{\operatorname{Int}}

\renewcommand{\setminus}{\smallsetminus}

\newcommand{\Area}{\operatorname{Area}}

\begin{document}

\title
%[ tropical geometry]
{Amoebas of algebraic varieties and tropical geometry}
\author{Grigory Mikhalkin}
\thanks{The author is partially supported by the NSF}
\address{Department of Mathematics\\
University of Utah\\
Salt Lake City, Utah 84112 USA}
\address{Department of Mathematics\\
University of Toronto\\
Toronto, Ontario M5S 3G3 Canada}
\address{St. Petersburg Branch of Steklov Mathematical Institute,
Fontanka 27, St. Petersburg, 191011 Russia}
\email{mikha@math.toronto.edu}

\maketitle

This survey consists of two parts.
Part 1 is devoted to amoebas. These are images of algebraic
subvarieties in $\C^n\supset\tor$ under the logarithmic
moment map. The amoebas have essentially piecewise-linear
shape if viewed {\em at large}. Furthermore, they degenerate
to certain piecewise-linear objects called {\em tropical varieties}
whose behavior is governed by algebraic geometry over the so-called
tropical semifield. Geometric aspects of tropical algebraic geometry
are the content of Part 2. We pay special attention to tropical curves.
Both parts also include relevant applications of the theories.
Part 1 of this survey is a revised and updated version of
the report \cite{M-survey}.

\part{AMOEBAS}
\section{Definition and basic properties of amoebas}
\subsection{Definitions}
Let $V\subset\tor$ be an algebraic variety. Recall that $\C^*=\C\setminus 0$ is the
group of complex numbers under multiplication.
Let $\Log:\tor\to\R^n$ be defined by $\Log(z_1,\dots,z_n)\to(\log|z_1|,\dots,\log|z_n|)$.
\begin{defn}[Gelfand-Kapranov-Zelevinski \cite{GKZ}]
The {\em amoeba} of $V$ is $\am=\Log(V)\subset\R^n$.
\end{defn}

\begin{figure}[h]
\label{amofline}
\centerline{
\psfig{figure=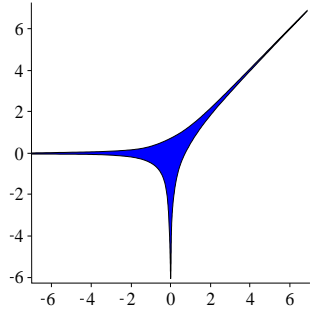,height=1.6in,width=1.6in}}
\caption{The amoeba of the line $\{ x+y+1=0\} \subset (\C^*)^2$.}
\end{figure}
%Figure \ref{amofline} is borrowed from \cite{MV}.

\begin{prop}[\cite{GKZ}]
The amoeba $\am\subset\R^n$ is a closed set with a non-empty complement.
\end{prop}

If $\C T\supset\tor$ is a closed $n$-dimensional toric variety and $\bar{V}\subset\C T$ is
a compactification of $V$ then we say that $\am$ is the amoeba of $\bar{V}$ (recall
that $\am$ is also the amoeba of $V=\bar{V}\cap\tor$).
Thus we can speak about amoebas of projective varieties once the coordinates
in $\cp^n$, or at least an action of $\tor$, is chosen.

If  $\C T$ is equipped with a $\tor$-invariant symplectic form then we can
also consider the corresponding moment map $\bar\mu:\C T\to\Delta$
(see \cite{At},\cite{GKZ}),
where $\Delta$ is the convex polyhedron associated to the toric variety
$\C T$ with the
given symplectic form. The polyhedron $\Delta$ is a subset of $\R^n$ but it is well defined
only up to a translation.
In this case we can also define the {\em compactified amoeba} of $\bar{V}$.
\begin{defn}[\cite{GKZ}]
The {\em compactified amoeba} of $V$ is $\bar\am=\bar\mu(V)\subset\Delta$.
\end{defn}

\begin{figure}[h]
\centerline{\psfig{figure=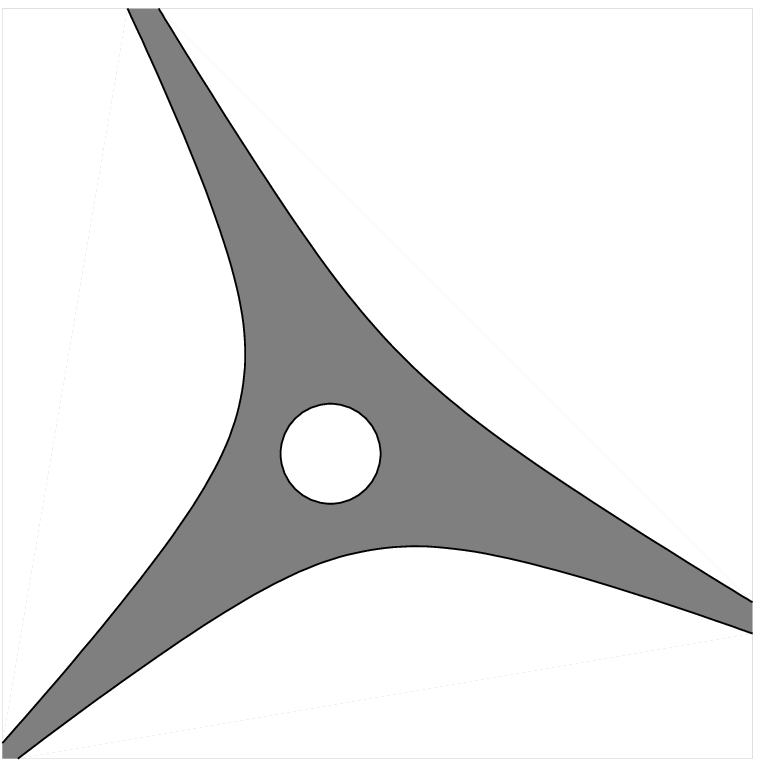,height=1in,width=1in}\hspace{1in}
\psfig{figure=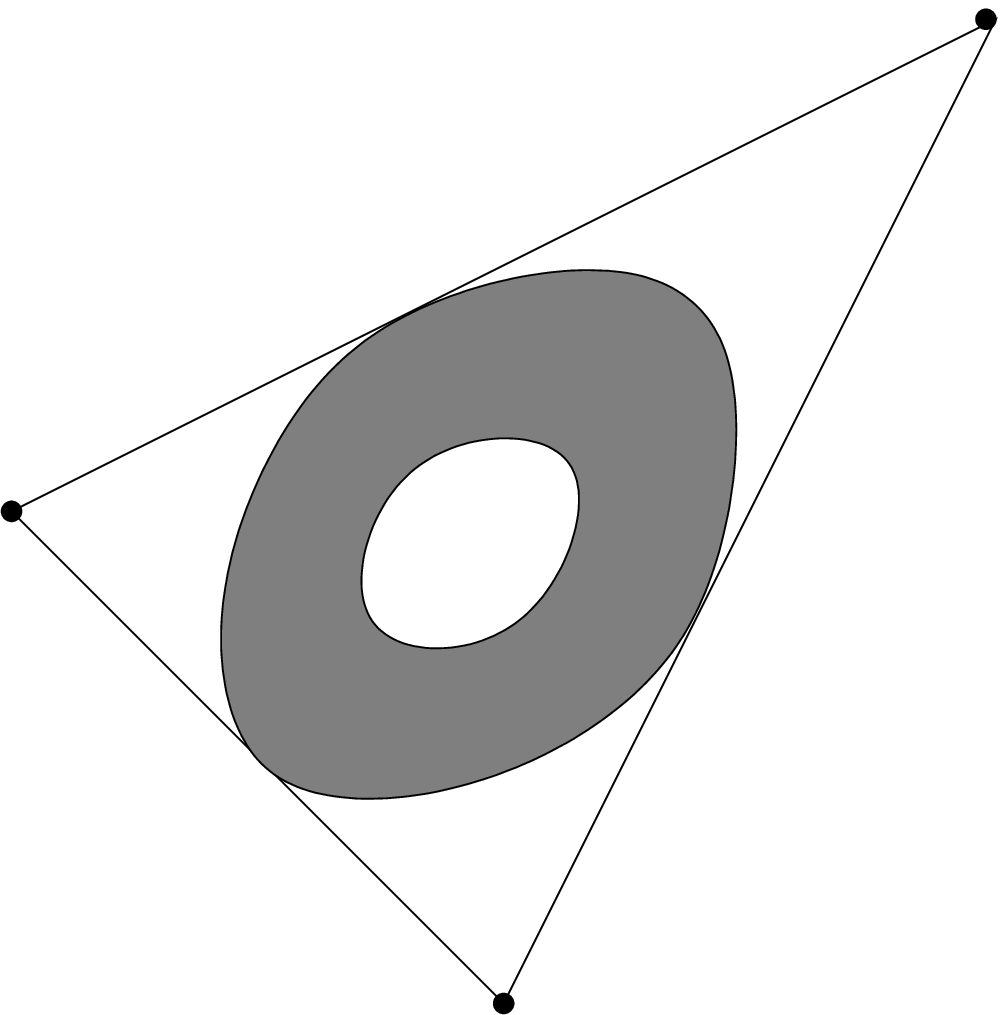,height=1in,width=1in}}
\caption{\cite{Mi} An amoeba $\am$ and a compactified amoeba $\bar{\am}$.}
\end{figure}
\begin{rem}
\label{repar}
Maps $\bar\mu|_{\tor}$ and $\Log$ are submersions and have the same real $n$-tori as fibers.
Thus $\am$ is mapped diffeomorphically onto $\bar\am\cap\Int\Delta$ under a
reparameterization of $\R^n$ onto $\Int\Delta$.
\end{rem}

%\subsection{Amoebas at infinity}
Using the compactified amoeba we can describe the behavior of $\am$
near infinity. Note that each face $\Delta'$ of $\Delta$ determines a
toric variety $\C T'\subset\C T$. Consider $\bar{V}'=\bar{V}\cap\C T'$.
Let $\bar\am'$ be the compactified amoeba of $\bar{V}'$.
\begin{prop}[\cite{GKZ}]
\label{natam}
We have
$\bar\am'=\bar\am\cap\Delta'$.
\end{prop}
This proposition can be used to describe the behavior of $\am\subset\R^n$ near infinity.

\subsection{Amoebas at infinity}
Consider a linear subspace $L\subset\R^n$ parallel to $\Delta'$ and
with $\dim L=\dim \Delta'$. Let $H\subset\R^n$ be a supporting hyperplane
for the convex polyhedron $\Delta$ at the face $\Delta'$, i.e. a hyperplane such
that $\Delta\cap H=\Delta'$.
Let $\stackrel{\to}{v}$ be an outwards normal vector to $H$.
Let $\am^{\Delta'}_t$, $t>0$, be the intersection of $L$ with the
result of translation of $\am$ by $-t\stackrel{\to}{v}$.

Recall that the {\em Hausdorff distance} between two closed sets
$A,B\subset\R^n$ is $$\dH(A,B)=\max\{\sup\limits{a\in A} d(a,B),
\sup\limits{b\in B} d(b,A)\},$$
where $d(a,B)$ is the Euclidean distance between a point $a$ and a set $B$.
We say that a sequence $A_t\subset\R^n$ converges to a set $A'$ when $t\to\infty$
{\em with respect to the Hausdorff metric on compacts in $\R^n$} if for any
compact $K\subset\R^n$ we have $\lim\limits_{t\to\infty}\dH(A_t\cap K,A'\cap K)=0$.

\begin{prop}
\label{asym}
The subsets $\am^{\Delta'}_t$ converge to $\am'$ when $t\to\infty$
with respect to the Hausdorff metric on compacts in $\R^n$.
\end{prop}

\begin{figure}[h]
\centerline{
\psfig{figure=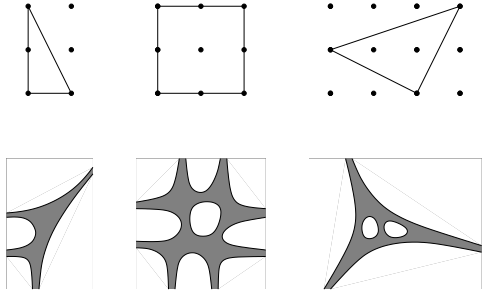,height=1.5in,width=2.5in}}
\caption{\label{3tent} Amoebas together with their Newton polyhedra.}
\end{figure}

This proposition can be informally restated in the case $n=2$ and $\dim V=1$.
In this case $\Delta$ is a polygon and the amoeba $\am$ develops ``tentacles"
perpendicular to the sides of $\Delta$ (see Figure \ref{3tent}). The number
of tentacles perpendicular to a side of $\Delta$ is bounded from
above by the integer length
of this side, i.e. one plus the number of the lattice points in the interior of the side.
%Intersections of $\am$ with affine subspaces of $\R^n$ which are sufficiently
%far from the origin look like thickenings of smaller-dimensional amoebas.

%Note that we may assume (by passing to a different toric variety $\C T$ if needed)
%that $V$ does not pass through the vertices of $\C T$, i.e. the fixed points of the $\tor$-action.
%Thus we get the following corollary.
\begin{coro}
For a generic choice of the slope of a line $\ell$ in $\R^n$ the intersection $\am\cap\ell$ is compact.
\end{coro}

%The advantage of considering the amoeba $\am$ over
%the compactified amoeba is in some geometric properties of $\am$.
%
%\subsection{Concavity of amoebas and homology of $\R^n\setminus\am$}
%\subsubsection{The case when $V\subset\tor$ is a hypersurface}
\subsection{Amoebas of hypersurfaces: concavity and topology of the complement}
\label{fpt}
Forsberg, Passare and Tsikh treated amoebas of hypersurfaces in \cite{FPT}.
In this case $V$ is a zero set of a single polynomial $f(z)=\sum\limits_ja_jz^j,
a_j\in\C$. Here we use the multiindex notations $z=(z_1,\dots,z_n)$,
$j=(j_1,\dots,j_n)\in\Z^n$ and $z^j=z_1^{j_1}\dots z_n^{j_n}$.
Let
\begin{equation}
\label{NP}
\Delta=\text{Convex hull} \{j\ |\ a_j\neq 0\}\subset\R^n
\end{equation}
be the Newton polyhedron of $f$.
\begin{thm}[Forsberg-Passare-Tsikh \cite{FPT}]
\label{thmfpt}
Each component of $\R^n\setminus\am$ is a convex domain in $\R^n$.
There exists a locally constant function
$$\ind:\R^n\setminus\am\to\Delta\cap\Z^n$$
which maps different components of the complement of $\am$ to
different lattice points of $\Delta$.
\end{thm}
\begin{cor}[\cite{FPT}]
\label{chislokomponent}
The number of components of $\R^n\setminus\am$ is never greater
then the number of lattice points of $\Delta$.
\end{cor}
Theorem \ref{thmfpt} and Proposition \ref{asym} indicate the dependence
of the amoeba on the Newton polyhedron.

%Figure \ref{am

The inequality of Corollary \ref{chislokomponent} is sharp.
This sharpness is a special case of Theorem \ref{komponenty}.
Also examples of amoebas with the maximal number
of the components of the complement are supplied by Theorem \ref{amlim}.

The concavity of $\am$ is equivalent to concavity of its boundary.
The boundary $\dd\am$ is contained in the critical value locus of $\Log|_V$.
The following proposition also takes care of some interior branches of this
locus.
\begin{prop}[\cite{Mi}]
\label{lc}
Let $D\subset\R^n$ be an open convex domain and $V'$ be a connected
component of $\Log^{-1}(D)\cap V$. Then $D\setminus\Log(V')$ is convex.
\end{prop}

\subsection{Amoebas in higher codimension: concavity}
The amoeba of a hypersurface is of full dimension in $\R^n$, $n>1$,
unless its Newton polyhedron $\Delta$ is contained in a line.
The boundary $\dd\am$ at its generic point is a smooth $(n-1)$-dimensional
submanifold. Its normal curvature form has no negative squares with
respect to the outwards normal (because
of convexity of components of $\R^n\setminus\am$). This property can be
generalized to the non-smooth points in the following way.
\begin{defn}
An open interval $D^1\subset L$, where $L$ is a straight line in $\R^n$,
is called a {\em supporting 1-cap} for $\am$ if
\begin{itemize}
\item $D^1\cap\am$ is non-empty and compact;
\item there exists a vector $\stackrel{\to}{v}\in\R^n$ such that the translation of $D^1$
by $\epsilon\hspace{-3pt}\stackrel{\to}{v}$ is disjoint from $\am$ for all sufficiently small $\epsilon>0$.
\end{itemize}
\end{defn}
The convexity of the components of $\R^n\setminus\am$ can be reformulated as stating
that {\em there are no 1-caps for $\am$}.

Similarly we may define higher-dimensional caps.
\begin{defn}
An open round disk $D^k\subset L$ of radius $\delta>0$ in a $k$-plane $L\subset\R^n$
is called a {\em supporting k-cap} for $\am$ if
\begin{itemize}
\item $D^k\cap\am$ is non-empty and compact;
\item there exists a vector $\stackrel{\to}{v}\in\R^n$ such that the translation of $D^k$
by $\epsilon\hspace{-3pt}\stackrel{\to}{v}$ is disjoint from $\am$ for all sufficiently small $\epsilon>0$.
\end{itemize}
\end{defn}

Consider now the general case, where $V\subset\tor$ is $l$-dimensional.
Let $k=n-l$ be the codimension of $V$.
The amoeba $\am$ is of full dimension in $\R^n$ if $2l\ge n$.
The boundary $\dd\am$ at its generic point is a smooth $(n-1)$-dimensional
submanifold. Its normal curvature form may not have more than
$k-1$ negative squares with respect to the outwards normal.
To see that note that a composition of $\Log|_V:V\to\R^n$ and
any linear projection $\R^n\to\R$ is a pluriharmonic function.

Note that this implies that there are no $k$-caps for $\am$
at its smooth points. It turns out that there are no $k$-caps
for $\am$ at the non-smooth points as well and also in the
case of $2l<n$ when $\am$ is $2l$-dimensional.
\begin{prop}[Local higher-dimensional concavity of $\am$]
\label{locconv}
If $V\subset\tor$ is of codimension $k$ then $\am$ does not
have supporting $k$-caps.
\end{prop}

A global formulation of convexity was
treated by Andr\'e Henriques \cite{Henr}.
\begin{defn}[Henriques \cite{Henr}]
A subset $\am\subset\R^n$ is called {\em $k$-convex} if for any $k$-plane $L\subset\R^n$
the induced homomorphism $H_{k-1}(L\setminus\am)\to H_{k-1}(\R^n\setminus\am)$
is injective.
\end{defn}
%\begin{thm}[Global higher-dimensional concavity of $\am$, cf. \cite{Henr}]
%\label{globconv}
%If $V\subset\tor$ is of codimension $k$ then $\am$ is $k$-convex.
%\end{thm}
Conjecturally the amoeba of a codimension $k$ variety in $\tor$ is $k$-convex.
A proof of a somewhat weaker version of this statement is contained in \cite{Henr}.
%Theorem \ref{globconv} can be deduced from its local version,
%Proposition \ref{locconv}.

\subsection{Amoebas in higher codimension: topology of the complement}
Recall that in the hypersurface case each component of $\R^n\setminus\am$
is connected and that there are not more than $\#(\Delta\cap\Z^n)$ such
components. The correspondence between the components of the complement
and the lattice points of $\Delta$ can be viewed as a cohomology class
$\alpha\in H^0(\R^n\setminus\am;\Z^n)$ whose evaluation on a point in each
component of $\R^n\setminus\am$ is the corresponding lattice point.

Similarly, when $V$ is of codimension $k$ there exists a natural class (cf. \cite{R2})
$$\alpha\in H^{k-1}(\R^n\setminus\am;H^k(T^n)),$$
where $T^n$ is the real $n$-torus, the fiber of $\Log$, $H^k(T^n)=H^k(\tor)$.
The value of $\alpha$ on each $(k-1)$-cycle $C$ in $\R^n\setminus\am$ and $k$-cycle $C'$
in $T^n$ is the linking number in $\C^n\supset\tor$ of $C\times C'$ and the closure of $V$.

The cohomology class $\alpha$ corresponds to the linking with the fundamental class
of $V$. Consider now the linking with smaller-dimensional homology of $V$.

Note that for an $l$-dimensional variety $V\subset\tor$ we have
 $H_j(V)=0$, $j>l$.
Similarly, $H^c_j(V)=0$, $j<l$, where $H^c$ stands for homology with closed support.
The linking number in $\R^n$ composed with $\Log:\tor\to\R^n$ defines the
following pairing
$$H^c_l(V)\times H_{k-1}(\R^n\setminus\am)\to\Z.$$
Together with the Poincar\'e duality between $H^c_l(V)$ and $H_l(V)$
this pairing defines the homomorphism
$$\iota: H_{k-1}(\R^n\setminus\am)\to H_l(V).$$
\begin{que}
\label{sparivanie}
Is $\iota$ injective?
\end{que}
Recall that a subspace $L\subset H_l(V)$ is called {\em isotropic} if the restriction
of the intersection form to $L$ is trivial.
\begin{prop}
\label{isotrop}
The image $\iota(H_{k-1}(\R^n\setminus\am))$ is isotropic in $H_l(V)$.
\end{prop}
\begin{rmk}
A positive answer to Question \ref{sparivanie} together with Proposition  \ref{isotrop}
would produce an upper bound for the dimension of $H_{k-1}(\R^n\setminus\am)$.
\end{rmk}
One may also define similar linking forms for $H_j(\R^n\setminus\am)$, $j\neq k-1$
(if $j>k-1$ then we can use ordinary homology $H_{n-j-1}(V)$
instead of homology with closed support) .

The answer to Question \ref{sparivanie} is currently unknown even in the
case when $V\subset(\C^*)^2$ is a curve. In this case $V$ is a Riemann surface
and it is defined by a single polynomial.
Let $\Delta$ be the Newton polygon of $V$.
The genus of $V$ is equal to the number of lattice points strictly inside $\Delta$
(see \cite{Kh}) while the number of punctures is equal to the number of lattice points
on the boundary of $\Delta$). Thus the dimension of a maximal isotropic subspace
of $H_1(V)$ is equal to $\#(\Delta\cap\Z^2)$ and Question \ref{sparivanie} agrees
with Corollary \ref{chislokomponent} for this case.

\section{Analytic treatment of amoebas}
\label{paru}
This section outlines the results obtained by Passare and Rullg{\aa}rd
in \cite{PR}, \cite{R1} and \cite{R2}.

We assume that $V\subset\tor$ is a hypersurface in this section.
Thus $V=\{f=0\}$ for a polynomial $f:\tor\to\C$ and we can consider
$\Delta\subset\R^n$, the Newton polyhedron of $V$ (see \ref{fpt}).

\subsection{The Ronkin function $N_f$}
Since $f$ is a holomorphic function,
$\log|f|:\tor\setminus V\to\R$ is a pluriharmonic function.
Furthermore, if we set $\log(0)=-\infty$ then
we have a plurisubharmonic function
$$\log|f|:\tor\to\R\cup\{-\infty\},$$
which is, obviously, strictly plurisubharmonic over $V$.
Recall that a function $F$ in a domain $\Omega\subset\C^n$ is called
plurisubharmonic if its restriction to any complex line $L$ is subharmonic,
i.e. the value of $F$ at each point $z\in L$ is smaller or equal than the
average of the value of $F$ along a small circle in $L$ around $z$.

Let $N_f:\R^n\to\R$ be the push-forward of $\log|f|$ under the map $\Log:\tor\to\R^n$,
i.e.
$$N_f(x_1,\dots,x_n)=\frac{1}{(2\pi i)^n}\int\limits_{\Log^{-1}(x_1,\dots,x_n)}\log|f(z_1,\dots,z_n)|
\frac{dz_1}{z_1}\wedge\dots\wedge\frac{dz_n}{z_n},$$
cf. \cite{Ro}.
This function was called {\em the Ronkin function} in \cite{PR}.
It is easy to see that it takes real (finite) values even over $\am=\Log(V)$ where
the integral is singular.

\begin{prop}[Ronkin-Passare-Rullg{\aa}rd \cite{PR}, \cite{Ro}]
\label{ronkin}
The function $N_f:\R^n\to\R$ is convex. It is strictly convex
over $\am$ and linear over each component of $\R^n\setminus\am$.
\end{prop}
This follows from plurisubharmonicity of $\log|f|:\tor\to\R$, its strict
plurisubharmonicity over $V$ and its pluriharmonicity in $\tor\setminus V$.
Indeed the convexity of a function in a connected real domain is just
a real counterpart of plurisubharmonicity. A harmonic function of one real variable has
to be linear and thus a function of several real variables is real-plurisubharmonic
if and only if it is convex. Over each connected component of $\R^n\setminus\am$
the function is
linear as the push-forward of a pluriharmonic function.

\begin{rem}
Note that just the existence of a convex function $N_f$, which is
strictly convex over $\am$ and linear over components of $\R^n\setminus\am$,
implies that each component of $\R^n\setminus\am$ is convex.
\end{rem}

Thus the gradient $\nabla N_f:\R^n\to\R^n$ is constant over each component $E$ of $\R^n\setminus\am$.
Recall the classical Jensen's formula in complex analysis
$$\frac{1}{2\pi i}\int\limits_{|z|=e^x}\log|f(z)|\frac{dz}{z}=Nx+\log|f(0)|-\sum\limits_{k=1}^{N}\log|a_k|,$$
where $a_1,\dots,a_N$ are the zeroes of $f$ in $|z|<e^x$, if $f(0)\neq 0$ and $f(z)\neq 0$ if $|z|=e^x$.
This formula implies that $\nabla N_f(E)\in\Z^n\cap\Delta$.
\begin{prop}[Passare-Rullg{\aa}rd \cite{PR}]
We have
$$\Int\Delta\subset\nabla N_f(\R^n)\subset\Delta,$$
where $\Int\Delta$ is the interior of the Newton polyhedron.
\end{prop}
Recall that Theorem \ref{thmfpt} associates a lattice point
to each component of $\R^n\setminus\am$.
\begin{prop}[\cite{PR}]
We have $$\nabla N_f(E)=\ind(E)$$
for each component $E$ of $\R^n\setminus\am$.
\end{prop}

\subsection{The spine of amoeba}

Passare and Rullg{\aa}rd \cite{PR} used $N_f$ to define {\em the spine} of amoeba.
Recall that $N_f$ is piecewise-linear on $\R^n\setminus\am$ and convex in $\R^n$.
Thus we may define a superscribed convex linear function $N^{\infty}_f$ by letting
$$N^{\infty}_f=\max_E N_E,$$
where $E$ runs over all components of $\R^n\setminus E$ and $N_E:\R^n\to\R$ is the linear
function obtained by extending $N_f|_E$ to $\R^n$ by linearity.

\begin{defn}[\cite{PR}]
\label{spine}
The spine $\s$ of amoeba is the corner locus of $N^{\infty}_f$, i.e.
the set of points in $\R^n$ where $N^{\infty}_f$ is not locally linear.
\end{defn}
Note that $\s\subset\am$ and that $s$ is a piecewise-linear polyhedral complex.
The following theorem shows that $\s$ is indeed a spine of $\am$ in the
topological sense.

\begin{figure}[h]
\label{am-spine}
\centerline{
\psfig{figure=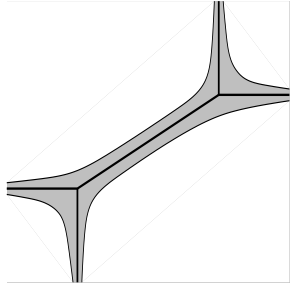,height=1.5in,width=1.5in}}
\caption{An amoeba and its spine.}
\end{figure}

\begin{thm}[\cite{PR}, \cite{R2}]
The spine $\s$ is a strong deformational retract of the amoeba $\am$.
\end{thm}
Thus each component of $\R^n\setminus\s$ (i.e. each maximal open domain
where $N^{\infty}_f$ is linear) contains a unique component of $\R^n\setminus\am$.
%The spine $\s$ can be considered as a certain non-Archimedean amoeba.

\subsection{Spine of amoebas and some functions on the space of complex polynomials}
Now we return to the study of the spine $\s\subset\am$ of a complex amoeba.
The spine $\s$ itself a certain amoeba over a non-Archimedean
field $K$. It does not matter what is the field $K$ as long as the corresponding hypersurface
over $K$ has the coefficients $a_j\in K$ with the correct valuations.
We can find these valuations from $N^{\infty}_f$ by taking its Legendre transform.
Since $N^{\infty}_f$ is obtained as a maximum of a finite number of linear function with
integer slopes its Legendre transform has a support on a convex lattice polyhedron $\Delta\subset\R^n$.
Let $c_{\alpha}\in\R$, $\alpha\in\Delta\cap\Z^n$ be the value of the Legendre transform of $N^{\infty}_f$ at $\alpha$.
To present $\s$ as a non-Archimedean amoeba we choose $a_j\in K$ such that $v(a_j)=c_{\alpha}$.

For each $\alpha\in\Delta\cap\Z^n$ let $U_{\alpha}$ be the space of all polynomials
whose Newton polyhedron is contained in $\Delta$ and whose amoeba contains a component
of the complement of index $\alpha$. The space of all polynomials whose Newton polyhedron
is contained in $\Delta$ is isomorphic to $\C^N$, where $N=\#(\Delta\cap\Z^n)$.
The subset $U_{\alpha}\subset\C^N$ is an open domain.
Note that $c_{\alpha}$ defines a real-valued function
on $U_{\alpha}$. This function was used by Rullg{\aa}rd \cite{R1}, \cite{R2} for
the study of geometry of $U_{\alpha}$.

\subsection{Geometry of $U_\alpha$}
Fix $\alpha\in\Delta\cap\Z^n$. Consider the following function in the space $\C^N$
of all polynomials $f$ whose Newton polyhedron is contained in $\Delta$
$$u_\alpha(f)=\inf\limits_{x\in\R^n}\frac{1}{(2\pi i)^n}
\int\limits_{\Log^{-1}(x)}\log|\frac{f(z)}{z^{\alpha}}|\frac{dz_1}{z_1}\wedge\dots\wedge\frac{dz_n}{z_n}, z\in\tor.$$
 Rullg{\aa}rd \cite{R1} observed that this function is plurisubharmonic in $\C^N$ while
pluriharmonic over $U_\alpha$. Indeed, over $U_\alpha$ there is a component
$E_\alpha\subset\R^n\setminus\am$ corresponding to $\alpha$ and $u_\alpha=\re\Phi_\alpha$,
where
$$\Phi_\alpha=\frac{1}{(2\pi i)^n}\int\limits_{\Log^{-1}(x)}\log(\frac{f(z)}{z^{\alpha}})
\frac{dz_1}{z_1}\wedge\dots\wedge\frac{dz_n}{z_n}, x\in E_\alpha$$
is a $(\C/2\pi i\Z)$-valued holomorphic function. Note that over $\Log^{-1}(E_\alpha)$
we can choose a holomorphic branch of $\log(\frac{f(z)}{z^{\alpha}})$ and that
$\Phi_\alpha$ does not depend on the choice of $x\in E_\alpha$.
Therefore, $U_\alpha$ is pseudo-convex.

Note that $U_\alpha$ is invariant
under the natural $\C^*$-action in $\C^N$. Let ${\mathcal C}\subset\cp^{N-1}$ be the complement
of the image of $U_\alpha$ under the projection $\C^N\to\cp^{N-1}$.
\begin{thm}[Rullg{\aa}rd \cite{R1}]
For any line $L\subset\cp^{n-1}$ the set $L\cap {\mathcal C}$ is non-empty and connected.
\end{thm}

The next theorem describes how the sets $U_\alpha$ with different $\alpha\in\Delta\cap\Z^n$
intersect. It turns out that for any choice of subdivision $\Delta\cap\Z^n=A\cup B$ with $A\cap B=\emptyset$
the sets $\bigcup\limits_{\alpha\in A} U_\alpha$ and $\C^N\setminus\bigcup\limits_{\beta\in B} U_\beta$
intersect. A stronger statement was found by Rullg{\aa}rd.
Let $A, B\subset\Delta\cap\Z^n$ be disjoint sets. The set $A\cup B\subset\Delta\cap\Z^n$ defines
a subspace $\C^{\#(A\cup B)}\subset\C^N$ .
\begin{thm}[\cite{R1}]
\label{komponenty}
For any $\#(A\cup B)$-dimensional space $L$ parallel to $\C^{\#(A\cup B)}$ the intersection
$L\cap\bigcup\limits_{\alpha\in A} U_\alpha\cap\C^N\setminus\bigcup\limits_{\beta\in B} U_\beta$
is non-empty.
\end{thm}

\subsection{The Monge-Amp\`ere measure and the symplectic volume}

\begin{defn}[Passare-Rullg{\aa}rd \cite{PR}]
\label{ma}
The Monge-Amp\`ere measure on $\am$ is the pull-back of the Lebesgue measure on $\Delta\subset\R^n$
under $\nabla N_f$.
\end{defn}
Indeed by Proposition \ref{ronkin} the Monge-Amp\`ere measure is well-defined. Furthermore,
we have the following proposition.

\begin{prop}[\cite{PR}]
The Monge-Amp\`ere measure has its support on $\am$.
The total Monge-Amp\`ere measure of $\am$ is $\Vol\Delta$.
\end{prop}

By Definition \ref{ma} the Monge-Amp\`ere measure is given by the determinant of the Hessian of $N_f$.
By convexity of $N_f$ its Hessian $\operatorname{Hess} N_f$
is a non-negatively defined matrix-valued function.
The trace of $\operatorname{Hess} N_f$ is the Laplacian of $N_f$,
it gives another natural measure supported on $\am$.
Note that $\omega=\sum\limits_{k=1}^n\frac{dz_k}{z_k}\wedge\frac{d\bar z_k}{\bar z_k}$
is a symplectic form on $\tor$ invariant with respect to the group structure.
The restriction $\omega|_V$ is a symplectic form on $V$. Its $(n-1)$-th power
divided by $(n-1)!$ is
a volume form called {\em the symplectic volume} on the $(n-1)$-manifold $V$.
\begin{thm}[\cite{PR}]
The measure on $\am$ defined by the Laplacian of $N_f$ coincides with the
push-forward of the symplectic volume on $V$, i.e. for any Borel set $A$
$$\int\limits_A\Delta N_f=\int\limits_{\Log^{-1}(A)\cap V}\omega^{n-1}.$$
\end{thm}

This theorem appears in \cite{PR} as a particular case of a computation for
{\em the mixed Monge-Amp\`ere operator}, the symmetric multilinear
operator associating a measure to $n$ functions $f_1,\dots,f_n$ (recall that
by our convention $n$ is the number of variables) and such
that its value on $f,\dots,f$ is the Monge-Amp\`ere measure from
Definition \ref{ma}. The total mixed Monge-Amp\`ere measure for $f_1,\dots,f_n$
is equal to the mixed volume of the Newton polyhedra of $f_1,\dots,f_n$
divided by $n!$.

Recall that this mixed volume divided by $n!$ appears in the Bernstein formula \cite{B}
which counts the number of common solutions of the system of equations $f_k=0$
(assuming that the corresponding hypersurfaces intersect transversely).
Passare and Rullg{\aa}rd found the following local analogue of the Bernstein formula
which also serves as a geometric interpretation of the mixed Monge-Amp\`ere measure.
Note that the complex torus $\tor$ acts on polynomials of $n$ variables. The value
of $t\in\tor$ on $f:\tor\to\C$ is the composition $f\circ t$ of the multiplication by $t$ followed by
application of $f$. In particular, the real torus $T^n=\Log^{-1}(0)\subset\tor$ acts on
polynomials of $n$ variables.
\begin{thm}[\cite{PR}]
\label{be}
The mixed Monge-Amp\`ere measure for $f_1,\dots,f_n$ of a Borel set $A\subset\R^n$
is equal to the average number of solutions of the system of equations $f_k\circ t_k=0$
in $\Log^{-1}(E)\subset\tor$, $t_k\in T^n$, $k=1,\dots,n$.
\end{thm}
The number of solution of this system of equations does not depend on $t_k$
as long as the choice of $t_k$ is generic. Thus Theorem \ref{be} produces the Bernstein
formula when $E=\R^n$.

\subsection{The area of a planar amoeba}
The computations of the previous subsection can be used to obtain an upper
bound on amoeba's area in the case when $V\subset(\C^*)^2$ is a curve.
With the help of Theorem \ref{be} Passare and Rullg{\aa}rd \cite{PR} showed
that in this case the Lebesgue measure on $\am$ is not greater than $\pi^2$
times the Monge-Amp\`ere measure. In particular we have the following theorem.
\begin{thm}[\cite{PR}]
\label{PR}
If $V\subset(\C^*)^2$ is an algebraic curve then
$$\operatorname{Area}\am\le\pi^2\operatorname{Area}\Delta.$$
\end{thm}
This theorem is specific for the case $\am\subset\R^2$.
Non-degenerate higher-dimensional amoebas of hypersurfaces have infinite volume.
This follows from Proposition \ref{asym} since the area of the cross-section at infinity
must be separated from zero.

\section{Some applications of amoebas}
\subsection{The first part of Hilbert's 16th problem}
\label{H16}
Most applications considered here are in the framework of Hilbert's 16th problem.
Consider the classical setup of its first part, see \cite{Hi}.
Let $\R\bar{V}\subset\rp^2$ be a smooth algebraic curve of degree $d$.
{\em What are the possible topological types of pairs $(\rp^2,\R\bar{V})$ for a given $d$?}

Since $\R \bar{V}$ is smooth it is homeomorphic to a disjoint union of circles.
All of these circles must be contractible in $\R P^2$ (such circles are called {\em the ovals})
if $d$ is even. If $d$ is odd then exactly one of these circles is non-contractible.
Therefore, the topological type of $(\rp^2,\R\bar{V})$ (also called {\em the topological
arrangement} of $\R\bar{V}$ in $\rp^2$) is determined by the number of components
of $\R\bar{V}$ together with the information on the mutual position of the ovals.

The possible number of components of $\R\bar{V}$ was determined by Harnack \cite{Ha}.
He proved that it cannot be greater than $\frac{(d-1)(d-2)}{2}+1$. Furthermore he proved
that for any number $$l\le\frac{(d-1)(d-2)}{2}+1$$ there exists a curve of degree $d$ with
exactly $l$ components as long as $l>0$ in the case of odd $d$ (recall that for odd $d$
we always have to have a non-contractible component).

Note that each oval separates $\rp^2$ into its {\em interior}, which is homeomorphic to a disk,
and its {\em exterior}, which is homeomorphic to a M\"obius band. If the interiors of the
ovals intersect then the ovals are called {\em nested}. Otherwise the ovals are called
{\em disjoint}. Hilbert's problem started from a question whether a curve of degree 6
which has 11 ovals (the maximal number according to Harnack) can have all of the ovals
disjoint. This question was answered negatively by Petrovsky \cite{P} who showed
that at least two ovals of a sextic must be nested if the total number of ovals is 11.

In general the number of topological  arrangements of curves of degree $d$
 grows exponentially with $d$. Even for small $d$ the number of the possible
types is enormous. Many powerful theorems restricting possible topological
arrangements were found for over 100 years of history of this problem,
see, in particular, \cite{P}, \cite{A}, \cite{R}, \cite{W}.
A powerful {\em patchworking}
construction technique \cite{V} counters these theorems.
The complete classifications is currently known for $d\le 7$, see \cite{V}.

The most restricted turn out to be curves with the maximal numbers of
components, i.e. with $l=\frac{(d-1)(d-2)}{2}+1$. Such curves were called
{\em M-curves} by Petrovsky. However, even for M-curves, the number of topological
arrangements grows exponentially with $d$.

The situation becomes different if we consider $\rp^2$ as a toric surface,
i.e. as a compactification of $(\R^*)^2$. Recall that $\rp^2\setminus (\R^*)^2$
consists of three lines $l_0$, $l_1$ and $l_2$ which can be viewed as
coordinate axes for homogeneous coordinates in $\rp^2$.
Thus we have three affine charts for $\rp^2$. The intersection of all
three charts is $(\R^*)^2\subset\rp^2$. We denote $\R V=\R\bar{V}\cap(\R^*)^2$.
The complexification $V\subset (\C^*)^2$ is the complex hypersurface defined
by the same equation as $\R V$. Thus we are in position to apply the
content of the previous sections of the paper to the amoeba of $V$.

In \cite{Mi} it was shown (with the help of amoebas) that for each $d$ the topological type of the
pair $(\rp^2,\R\bar{V})$ is unique as long as the curve $\R\bar{V}$ is
maximal in each of the three affine charts of $\rp^2$. Furthermore,
the diffeomorphism type of the triad $(\rp^2;\R\bar{V},l_0\cup l_1\cup l_2)$
is unique. In subsection
\ref{krivye} we formulate this maximality condition and sketch the
proof of uniqueness. A similar statement holds for curves in other
toric surfaces. The Newton polygon $\Delta$ plays then the r\^ole of
the degree $d$.
%In subsections \ref{poverxnosti} and \ref{hd} we describe an analogous
%but weaker statement towards uniqueness.

\subsection{Relation to amoebas: the real part $\R V$ as a subset of the critical locus of $\Log|_V$
and the logarithmic Gauss map}

Suppose that the hypersurface $V\subset\tor$ is defined over real numbers (i.e. by a polynomial
with real coefficients).
Denote its real part via $\R V=V\cap\rtor$.
We also assume that $V$ is non-singular.
Let $F\subset V$ be the critical locus of the map $\Log|_V:V\to\R^n$.
It turns out that the real part $\R V$ is always contained in $F$.
\begin{prop}[Mikhalkin \cite{Mi}]
\label{RVF}
$\R V\subset F$.
\end{prop}
This proposition indicates that the amoeba must carry some information about $\R V$.
%Hypersurfaces with $\R V=F$ are maximal in a certain sense and have a peculiar topology.
The proof of this proposition makes use of the {\em logarithmic Gauss map}.

Note that since $\tor$ is a Lie group there is a canonical trivialization
of its tangent bundle. If $z\in\tor$ then the multiplication by $z^{-1}$
induces an isomorphism $T_z\tor\approx T_1\tor$ of the tangent bundles at $z$
and $1=(1,\dots,1)\in\tor$.
\begin{defn}[Kapranov \cite{KaGauss}]
The {\em logarithmic Gauss map} is a map
$$\gamma:V\to\cp^{n-1}.$$ It sends each point $z\in V$ to the image
of the hyperplane $T_zV\subset T_z\tor$  under the canonical
isomorphism $T_z\tor\approx T_1\tor=\C^n$.

The map $\gamma$ is a composition of a branch of a holomorphic logarithm
$\tor\to\C^n$ defined locally up to translation by $2\pi i$ with the usual Gauss
map of the image of $V$.
We may define $\gamma$ explicitly in terms of the defining
polynomial $f$ for $V$ by logarithmic differentiation formula.
If $z=(z_1,\dots,z_n)\in V$ then
$$\gamma(z)=[<\nabla f,z>]=[\frac{\dd f}{\dd z_1}z_1:\dots:\frac{\dd f}{\dd z_n}z_n]\in\cp^{n-1}.$$
\end{defn}

\begin{lem}[\cite{Mi}]
$F=\gamma^{-1}(\rp^{n-1})$
\end{lem}
To justify this lemma we recall that $\Log:\tor\to\R^n$ is a smooth fibration
and $V$ is non-singular.
Thus $z\in V$ is critical for $\Log|_V$ if and only if the tangent vector
space to $V$ and the tangent vector space to the fiber torus
$\gamma^{-1}(\gamma(z))$ intersect along an $(n-1)$-dimensional
subspace. Such points are mapped to real points of $\cp^{n-1}$ by $\gamma$.

Note that this lemma implies Proposition \ref{RVF}. If $V$ is defined over
$\R$ then $\gamma$ is equivariant with respect to the complex conjugation
and maps $\R V$ to $\rp^{n-1}$.

\subsection{Compactification: a toric variety associated to a hypersurface in $\tor$}
A hypersurface $V\subset\tor$ is defined by a polynomial $f:\C^n\to C$.
If the coefficients of $f$ are real then we define the real part of $V$
by $\R V=V\cap\rtor$. Recall that the Newton polyhedron $\Delta\subset\R^n$
of $V$ is an integer convex polyhedron obtained as
the convex hull of the indices of monomials participating in $f$,
see \eqref{NP} in subsection \ref{fpt}.

Let $\C T_\Delta\supset\tor$ be the toric variety corresponding to $\Delta$,
see e.g. \cite{GKZ} and let $\R T_\Delta\supset\rtor$ be its real part.
We define $\bar{V}\subset\C T_\Delta$ as the closure of $V$ in $\C T_\Delta$
and we denote via $\R\bar{V}$ its real part.

Note that $\bar{V}$ may be singular even if $V$ is not. Nevertheless $\C T_\Delta$
is, in some sense, the best toric compactification of $\tor$ for $V$. Namely, $\bar{V}$ does not
pass via the points of $\C T_\Delta$ corresponding to the vertices of $\Delta$ and
therefore it does not have singularities there. Furthermore, $\C T_\Delta$ is minimal
among such toric varieties, since $\bar{V}$ intersect any line in $\C T_\Delta$
corresponding to an edge of $\Delta$.

Thus we may naturally compactify the pair $(\tor,V)$ to the pair $(\C T_\Delta,\bar{V})$.
In such a setup the polyhedron $\Delta$ plays the r\^ole of the degree in $\C T_\Delta$.
Indeed, two integer polyhedra $\Delta$ define the same toric variety $\C T_\Delta$
if their corresponding faces are parallel. But the choice of $\Delta$ also fixes
the homology class of $\bar{V}$ in $H_{2n-2}(\C T_\Delta)$.

The simplest example is the projective space $\cp^n$. The corresponding $\Delta$ is,
up to translation and the action of $SL_n(\Z)$ the simplex defined by equations
$z_j>0$, $z_1+\dots+z_n<d$. Thus in this case $\Delta$ is parameterized by a single
natural number $d$ which is the degree of $\bar{V}\subset\cp^n$.

\subsection{Maximality condition for $\R V$}
%\subsubsection{Smith inequality}
The inequality $l\le\frac{(d-1)(d-2)}{2}$ discovered by Harnack for the number $l$
of components of a curve $\R\bar{V}$ is a part of a more general {\em Harnack-Smith
inequality}.
Let $X$ be a topological space and let $Y$ be the fixed point set of a of a continuous
involution on $X$. Denote by $b_*(X;\Z_2)=\dim H_*(X;\Z_2)$ the total $\Z_2$-Betti number of $X$.
\begin{thm}[P. A. Smith, see e.g. the appendix in \cite{W}]
\label{smith}
$$b_*(Y;\Z_2)\le b_*(X;\Z_2).$$
\end{thm}
\begin{cor}
$b_*(\R\bar{V};\Z_2)\le b_*(\bar{V};\Z_2)$, $b_*(\R{V};\Z_2)\le b_*({V};\Z_2).$
\end{cor}
Note that Theorem \ref{smith} can also be applied to pairs which consist of
a real variety and real subvariety and other similar objects.
\begin{defn}[Rokhlin \cite{R}]
A variety $\R\bar{V}$ is called an {\em M-variety} if $$b_*(\R\bar{V};\Z_2)= b_*(\bar{V};\Z_2).$$
\end{defn}
E.g. if $\bar{V}\subset\cp^2$ is a smooth curve of degree $d$ then $\bar{V}$ is a Riemann surface
of genus $g=\frac{(d-1)(d-2)}{2}$.
Thus $b_*(\bar{V};\Z_2)=2+2g$. On the other hand, $b_*(\R\bar{V};\Z_2)=2l$, where $l$ is the number
of (circle) components of $\R\bar{V}$.
%Letting $X=\bar{V}$ and $Y=\R \bar{V}$ in Theorem \ref{smith} we get $l\le g+1=\frac{(d-1)(d-2)}{2}+1$.

%\subsubsection{Maximality condition}
Let $\R V\subset\rtor$ be an algebraic hypersurface, $\Delta$ be its Newton polyhedron,
$\R T_{\Delta}$ be the toric variety corresponding to $\Delta$ and $\R\bar{V}\subset\R T_\Delta$
the closure of $\R V$ in $\R T_\Delta$.
We denote with $V\subset\tor$ and $\bar{V}\subset\C T_\Delta$ the complexifications of these
objects.
Recall (see e.g. \cite{GKZ}) that each (closed) $k$-dimensional face $\Delta'$ of $\Delta$ corresponds to
a closed $k$-dimensional toric variety $\R T_{\Delta'}\subset\R T_\Delta$ (and, similarly,
$\C T_{\Delta'}\subset\C T_\Delta$). The intersection $V_{\Delta'}=\bar{V}\cap\C T_{\Delta'}$ is itself
a hypersurface in the $k$-dimensional toric variety $\C T_{\Delta'}$ with the Newton polyhedron $\Delta'$.
Its real part is $\R V_{\Delta'}=V_{\Delta'}\cap\R\bar{V}$.

Denote with $\operatorname{St}\Delta'\subset\dd\Delta$ the union of all the closed faces of $\Delta$
containing $\Delta'$. Denote
$V_{\operatorname{St}\Delta'}=\bigcup\limits_{\Delta''\subset\operatorname{St}\Delta'}V_{\Delta'}$
and $\R V_{\operatorname{St}\Delta'}=V_{\operatorname{St}\Delta'}\cap\R T_\Delta$.

\begin{defn}
%[Mikhalkin \cite{M-r}]
\label{max}
A hypersurface $\R\bar{V}\subset\C T_{\Delta}$ is called {\em torically maximal} if the following conditions
hold
\begin{itemize}
\item $\R\bar{V}$ is an M-variety, i.e. $b_*(\R\bar{V};\Z_2)=b_*(\bar{V};\Z_2)$;
\item the hypersurface $\bar{V}\cap\C T_{\Delta'}\subset\C T_{\Delta}$ is torically maximal for each
face $\Delta'\subset\Delta$ (inductively we assume that this notion is already defined in smaller
dimensions);
\item for each face $\Delta'\subset\Delta$ we have
$b_*(\R V\cup\R V_{\operatorname{St}\Delta'},\R V_{\operatorname{St}\Delta'};\Z_2)=
b_*(V\cup V_{\operatorname{St}\Delta'},V_{\operatorname{St}\Delta'};\Z_2)$.
\end{itemize}
\end{defn}

Consider a linear function $h:\R ^n\to\R$. A facet $\Delta'\subset\Delta$ is called {\em negative}
with respect to $h$ if the image of its outward normal vector under $h$ is negative.
We define
$\C T^-=\bigcup\limits_{\text{negative}\ \Delta'}\C T_{\Delta'}.$
In these formula we take the union over all the closed facets $\Delta'$ negative
with respect to $h$.
Let $V^-=\bar{V}\cap\C T^-$ and $\R V^-=V^-\cap\R\bar{V}$.

We call a linear function $h:\R^n\to\R$ generic if its kernel does not contain vectors orthogonal
to facets of $\Delta$.
\begin{prop}
If a hypersurface $\R\bar{V}\subset\R T_\Delta$ is torically maximal then for any generic linear function $h$
we have
$$b_*(\R V\cup\R V^-,\R V^-;\Z_2)=b_*(V\cup V^-,V^-;\Z_2).$$
\end{prop}

\subsection{Curves in the plane}
\label{krivye}
\subsubsection{Curves in $\rp^2$ and their bases}
%Let $\R V\subset (\R^*)^2$ be a curve with the Newton polygon $\Delta$.
%Consider the toric surface $\C T_\Delta$ corresponding to $\Delta$.
%We say that $\R V$ is generic if $\R\bar{V}\subset\R T_\Delta$ is non-singular
%and intersects transversely every
Note that if $\R V\subset (\R^*)^2$ is a torically maximal curve then
the number of components of $\R\bar{V}$ coincides with the genus of $\C \bar{V}$.
In other words (cf. \ref{H16}) $\R\bar{V}$ is an M-curve.

We start by reformulating the maximality condition of Definition \ref{max} for
the case of curves in the projective plane.
Let $\R C\subset\rp^2$ be a non-singular curve of degree $d$.

\begin{figure}[h]
\label{bbasa}
\centerline{
\psfig{figure=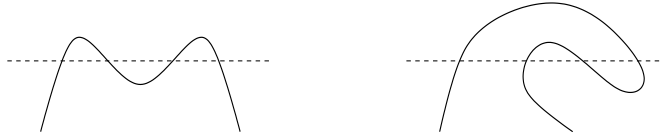,height=0.7in,width=3.5in}}
\caption{Possible bases for a real quartic curve.}
\end{figure}

\begin{defn}[Brusotti \cite{Br}]
Let $\alpha$ be an arc (i.e. an embedded closed interval) in $\R C$.
The arc $\alpha$ is called a {\em base} (or a {\em base of rank 1}, see \cite{Br})
if there exists a line $L\subset\rp^2$ such that the intersection $L\cap\alpha$
consists of $d$ distinct points.
\end{defn}

Note if three lines $L_1,L_2,L_3$ in $\rp^2$ are generic, i.e. they do not pass through
the same point, then =$\rp^2\setminus (L_1\cup L_2\cup L_3)=(\R^*)^2$.
We call such $(\R^*)^2$ a {\em toric chart} of $\rp^2$.
Thus $\R V=\R C\setminus (L_1\cup L_2\cup L_3)$ is a curve in $(\R^*)^2$.
If $\R C$ does not pass via $L_j\cap L_k$ then the Newton polygon of $\R V$
(for any choice of coordinates $(x,y)$ in $(\R^*)^2$ extendable to affine coordinates
in $\R^2=\rp^2\setminus L_j$ for some $j$) is the triangle $\Delta_d=\{x\ge 0\}\cap
\{y\ge 0\}\cap\{x+y\le d\}$.
\begin{prop}[Mikhalkin \cite{M-r}]
The curve $\R C\subset\rp^2$ is maximal in some toric chart of $\rp^2$
if and only if $\R C$ is an M-curve with three disjoint bases.
\end{prop}
Many M-curves with one or two disjoint bases are known (see e.g. \cite{Br}).
%In fact $\cite{Br}$ implies that the number of M-curves of degree
However there is (topologically) only one known example of curve with three disjoint bases,
namely the first M-curve constructed by Harnack \cite{Ha}. Theorem \ref{urp2}
asserts that this example is the only possible.

\begin{defn}[simple Harnack curve in $\rp^2$, cf. \cite{Ha}, \cite{MR}]
\label{har}
A non-singular curve $\R C\subset\rp^2$ of degree $d$ is called a (smooth) simple Harnack curve if it is an M-curve and
\begin{itemize}
\item all ovals of $\R C$ are disjoint (i.e. have disjoint interiors, see \ref{H16}) if $d=2k-1$ is odd;
\item one oval of $\R C$ contains $\frac{(k-1)(k-2)}{2}$ ovals in its interior while all other ovals
are disjoint if $d=2k$ is even.
\end{itemize}
\end{defn}

\begin{figure}[h]
\label{deg10}
\centerline{
\psfig{figure=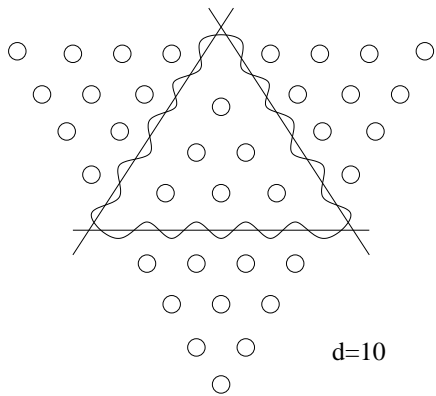,height=2.6in,width=2.6in}}
\caption{\cite{Mi} A simple Harnack curve.}
\end{figure}

\begin{thm}[\cite{Mi}]
\label{urp2}
Any smooth M-curve $\R C\subset\rp^2$ with at least three base is a simple Harnack curve.
\end{thm}
There are several topological arrangements of M-curves with fewer than 3 bases
for each $d$ (in fact, their number grows exponentially with $d$). There is a unique
(Harnack) topological arrangement of an M-curve with 3 bases by Theorem \ref{urp2}.
In the same time 3 is the highest number of bases an M-curve of sufficiently high degree
can have as the next theorem shows.
\begin{thm}[\cite{Mi}]
No M-curve in $\rp^2$ can have more than 3 bases if $d\ge 3$.
\end{thm}

\subsubsection{Curves in real toric surfaces}
Theorem \ref{urp2} has a generalization applicable to other toric surfaces.
%\begin{defn}[simple Harnack curve in $(\R^*)^2$, cf. \cite{Ha}, \cite{MR}]
%A non-singular curve $\R C\subset(\R^*)^2$ of degree $d$ is called a (smooth) simple Harnack curve
%if it is torically maximal.
%\end{defn}
Let $\R V\subset(\R^*)^2$ be a curve with the Newton polygon $\Delta$. The sides of $\Delta$ correspond
to lines $L_1,\dots,L_n$ in $\R T_\Delta$. We have $\R V=\R\bar{V}\setminus(L_1\cup\dots\cup L_n$.
\begin{thm}[\cite{Mi}]
\label{kriv}
The topological arrangement of a torically maximal curve is unique for each $\Delta$.
More precisely, the topological type of the triad $(\R T_\Delta;\R\bar{V},L_1\cup\dots\cup L_n)$
and, in particular, the topological type of the pair $((\R^*)^2,\R V)$ depends only on $\Delta$
as long as $\R V$ is a torically maximal curve.
\end{thm}

A torically maximal curve $\R\bar{V}$ is a counterpart of a simple Harnack curve for $\R T_\Delta$.
All of its components except for one are ovals with disjoint interiors. The remaining component
is not homologous to zero unless $\Delta$ is even (i.e. obtained from another lattice polygon
by a homotethy with coefficient 2). If $\Delta$ is even the remaining component is also an oval
whose interior contains $g(V)$ ovals of $\R V$. Recall that, by Khovanskii's formula \cite{Kh},
$g(V)$ coincides with the number of lattice points in the interior of $\Delta$.

\begin{thm}[Harnack, Itenberg-Viro \cite{Ha}, \cite{IV}]
For any $\Delta$ there exists a curve $\R V\subset(\R^*)^2$ which is torically maximal
and has $\Delta$ as its Newton polygon.
\end{thm}
As in Definition \ref{har} we call such curves {\em simple Harnack curves}, cf. \cite{MR}.

\subsubsection{Geometric properties of algebraic curves in $(\R^*)^2$}
It turns out that the simple Harnack curves have peculiar geometric properties,
but they are better seen after a logarithmic reparameterization $\Log|_{(\R^*)^2}:(\R^*)^2\to\R^2$.
A point of $\R V$ is called a logarithmic inflection point if it corresponds to an inflection
point of $\Log(\R V)\subset\R^2$ under $\Log$.

\begin{thm}[\cite{Mi}]
\label{logkriv}
The following conditions are equivalent.
\begin{itemize}
\item $\R V\subset(\R^*)^2$ is a simple Harnack curve.
\item $\R V\subset(\R^*)^2$ has no real logarithmic inflection points.
\end{itemize}
\end{thm}

\begin{rem}
Recall that by Proposition \ref{RVF} $\Log(\R V)$ is contained in the critical value
locus of $\Log|_V$. The map $\Log|_V:V\to\R^2$ is a surface-to-surface map in
our case and its most generic singularities are folds. By Proposition \ref{lc} the folds
are convex. Thus a logarithmic inflection point of $\R V$ must correspond to a
higher singularity of $\Log|_V$.

In \cite{Mi} it was stated that there are two types of stable (surviving small
deformations of $\R V$) logarithmic inflection points of $\R V$.
Here we'd like to correct this statement.
Only one of these two types is genuinely stable.
The first type (see Figure \ref{junct}), called
{\em junction}, corresponds to an intersection of $\R V$ with a
branch of imaginary folding curve. A junction logarithmic
inflection point can be found at the curve $y=(x-1)^2+1$. Note that the image of
the imaginary
folding curve under the complex conjugation is also a folding curve.
Thus over its image we have a double fold.

The second type, called {\em pinching}, corresponds to intersection of $\R V$ with
a circle $E\subset V$ that gets contracted by $\Log$.
%Such circles $E$ survive
%if we deform $V$ in the class of hypersurfaces with real coefficients but disappear
%under a generic small perturbation if we allow the coefficients to become imaginary.
The circle $E$ intersect $\R V$ at two points. These points belong to different
quadrants of $(\R^*)^2$, but have the same absolute values of their coordinates.
Both of these points are logarithmic inflection points.
The pinching is not stable even in the class of real deformations.
A small perturbation breaks it to two junctions with a corner of
two branches of the amoeba as in Figure \ref{pinching}.

\begin{figure}[h]
\centerline{
\psfig{figure=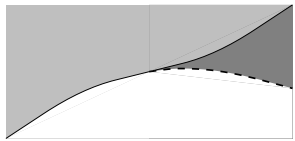,height=.7in,width=1.5in}}
\caption{\label{junct} A junction point.}
\end{figure}

\begin{figure}[h]
\centerline{
\psfig{figure=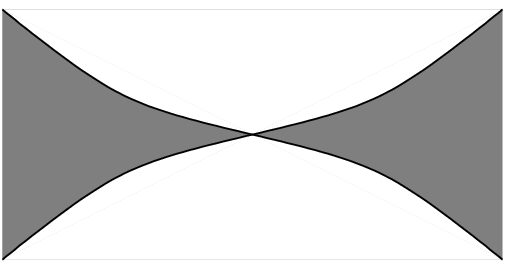,height=.7in,width=1.4in}\hspace{1in}
\psfig{figure=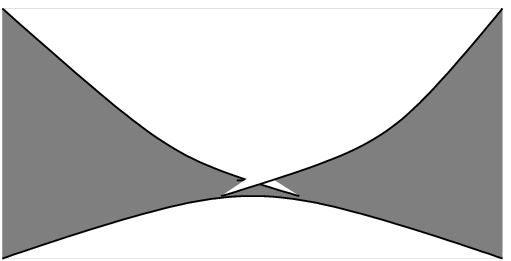,height=.7in,width=1.4in}}
\caption{\label{pinching} Deformation of a pinching point into two junction points.}
\end{figure}
\end{rem}

%Note that the boundary $\dd\am$ of amoeba must
\begin{prop}
The logarithmic image $\Log(\R V)$ is trivial in the closed support homology
group $H^c_1(\R^2)$.
\end{prop}
Thus the curve $\Log(\R V)$ spans a surface in $(\R^*)^2$.
%Thus there exists a closed domain ${\mathcal B}\subset\R^2$
%such that ${\mathcal B}\supset\Log(\R V)$ and $\dd{\mathcal B}\subset\Log(\R V)$.
%Theorem \ref{PR} has the following corollary.
Theorem \ref{PR} has the following corollary.
\begin{cor}
The area of any region spanned by branches of $\Log(\R V)$ is smaller
than $\Area\Delta$.
\end{cor}

The situation is especially simple for
the logarithmic image of a simple Harnack curve.
\begin{prop}[\cite{Mi}]
\label{shemb}
If $\R V$ is a simple Harnack curve then $\Log|_{\R V}$ is an embedding
and $\Log\R V=\dd\am$.
\end{prop}

Thus in this case $\am$ coincides with the region spanned by the whole curve $\Log(\R V)$.
Furthermore, in \cite{MR} it was shown that simple Harnack curves maximize the
area of this region.
\begin{thm}[Mikhalkin-Rullg{\aa}rd, \cite{MR}]
If $\R V$ is a simple Harnack curve then
$\Area\am=\Area\Delta$.
\end{thm}

In the opposite direction we have the following theorem.
We say that a curve $V\subset(\C^*)^2$ is real up to translation if there exists
$a\in(\C^*)^2$ such that $aV$ is defined by a polynomial with real coefficients.
We denote the corresponding real part with $\R V$.
(Note that in general this real part might depend on the choice of translation.)
\begin{thm}[\cite{MR}]
If $\Area\am=\Area\Delta>0$ and $V$ is non-singular and transverse to the
lines (coordinate axes) in $\C T_\Delta$ corresponding to the sides of $\Delta$
then $V$ is real up to translation in a unique way and $\R V$ is a simple Harnack curve.
\end{thm}
Furthermore, in \cite{MR} it was shown that the only singularities that $V$ can have
in the case $\Area\am=\Area\Delta>0$ are ordinary real isolated double points.

\subsection{A higher-dimensional case}
\label{poverxnosti}
\subsubsection{Surfaces in $({\mathbb R}^*)^3$}
Let $\R V\subset(\R^*)^3$ be an algebraic surface with the Newton
polyhedron $\Delta\subset\R^3$. Let $\R\bar{V}\subset\R T_\Delta$
be its compactification.

%\begin{prop}
Recall (see Definition \ref{max}) that if $\R V$ is a torically maximal surface then
$b_*(\R\bar{V};\Z_2)=b_*(\bar{V};\Z_2)$, i.e. $\R\bar{V}$ is an M-surface.
%\end{prop}

\begin{thm}[\cite{M-r}]
\label{pov}
Given a Newton polyhedron $\Delta$ the topological type of
a torically maximal surface $\R\bar{V}\subset\R T_\Delta$ is unique.
\end{thm}
To describe the topological type of $\R\bar{V}$ it is useful to  compute
the total Betti number $b_*(\bar{V};\Z_2)$ in terms of $\Delta$. Note that
by the Lefschetz hyperplane theorem $b_*(\bar{V};\Z_2)=\chi(\bar{V})$.

We denote by $\Area\dd\Delta$ the total area of the faces of $\Delta$.
Each of these faces sits in a plane $P\subset\R^3$. The intersection  $P\cap\Z^3$
determines the area form on $P$. This area form is
translation invariant and such that the area of the smallest lattice
parallelogram is 1.

Similarly we denote by $\operatorname{Length}\operatorname{Sk}^1\Delta$
the total length of all the edges of $\Delta$. Again, each edge sits in a line $L\subset\R^3$.
The intersection $L\cap\Z^3$ determines the length on $L$ by setting the length of
the smallest lattice interval 1.

\begin{prop}
$b_*(V;\Z_2)=6\Vol\Delta-2\Area\dd\Delta+\operatorname{Length}\operatorname{Sk}^1\Delta.$
\end{prop}
This proposition follows from Khovanskii's formula \cite{Kh}.

\begin{thm}[\cite{M-r}]
\label{ds}
A torically maximal surface $\R\bar{V}$ consists of $p_g+1$ components, where $p_g$
is the number of points in the interior of $\Delta$. There are $p_g$
components homeomorphic to 2-spheres and contained in $(\R^*)^3$.\
These spheres bound disjoint spheres in $(\R^*)^3$.
The remaining component is homeomorphic to
\begin{itemize}
\item a sphere with
$b_*(V;\Z_2)-2p_g(V)-2$ M\"obius bands in the case when $\Delta$ is odd
(i.e. cannot be presented as $2\Delta'$ for some lattice polyhedron $\Delta'$);
\item a sphere with $\frac12 b_*(V;\Z_2)-p_g(V)-1$ handles in the case $\Delta$ is even.
\end{itemize}
\end{thm}

\begin{rem}
\label{bertrand}
Not for every Newton polyhedron $\Delta$ a torically maximal surface $\R V\subset(\R^*)^3$
exists. The following example is due to B. Bertrand. Let $\Delta\subset\R^3$ be the convex hull
of $(1,0,0)$, $(0,1,0)$, $(1,1,0)$ and $(0,0,2k+1)$. If $k>0$ then there is no M-surface $\R\bar{V}$
with the Newton polyhedron $\Delta$. In particular, there is no torically maximal surface $\R V$
for $\Delta$.
\end{rem}

\begin{exa}
There are 3 different topological types of smooth M-quartics in $\rp^3$ (see \cite{Kha}).
They realize all topological possibilities for maximal real structures on abstract
K3-surfaces. Namely, such real surface may be homeomorphic to
\begin{itemize}
\item the disjoint union of 9 spheres and a surface of genus 2;
\item the disjoint union of 5 spheres and a surface of genus 6;
\item the disjoint union of a sphere and a surface of genus 10.
\end{itemize}
Theorem \ref{ds} asserts that only the last type can be a torically maximal quartic in $\rp^3$.
More generally, only the last type can be a torically maximal surface is a toric 3-fold
$\R T_\Delta$.
\end{exa}

%\ignore{
\subsubsection{Geometric properties of maximal algebraic surfaces in $(\R^*)^3$}
Recall the classical geometric terminology.
Let $S\subset\R^3$ be a smooth surface.
We call a point $x\in S$ {\em elliptic}, {\em hyperbolic} or {\em parabolic}
if the Gauss curvature of $S$ at $x$ is positive, negative or zero.
\begin{rem}
Of course we do not actually need to use the Riemannian metric on $S$
do define these points. Here is an equivalent definition without referring
to the curvature. Locally near $x$ we can present $S$ as the graph of a function
$\R^2\to\R$. If the Hessian form of this function at $x$ is degenerate then we
call $x$ parabolic. If not, the intersection of $S$ with the tangent plane at $x$
is a real curve with an ordinary double point in $x$. If this point is isolated
we call $x$ elliptic. If it is an intersection of two real branches of the curve
we call it hyperbolic.
\end{rem}

We say that a point $x\in\R V\subset (\R^*)^3$ is {\em logarithmically} elliptic,
hyperbolic or parabolic if it maps to such point under
$\Log|_{(\R^*)^3}:(\R^*)^3\to\R^3$.

Generically for a smooth surface in $\R^3$ the parabolic locus, i.e. the set of parabolic points,
is a 1-dimensional curve. So is the logarithmic parabolic locus for a surface in $(\R^*)^3$.
In a contrast to this we have the following theorem for torically maximal surfaces.
Note that torically maximal surfaces form an open subset in the space of all surfaces
with a given Newton polyhedron.

\begin{thm}[\cite{M-r}]
\label{parpov}
The logarithmic parabolic locus of a torically maximal surface
consists of a finite number of points.
\end{thm}
Note that such a zero-dimensional locus cannot separate the surface $\R V$.
Thus each component of $\R V$ is either logarithmically elliptic (all its points except finitely
many are logarithmically elliptic) or logarithmically hyperbolic
(all its points except finitely many are logarithmically hyperbolic).
\begin{cor}[\cite{M-r}]
Every compact component of $\R V$ is diffeomorphic to a sphere.
\end{cor}
This corollary is a part of Theorem \ref{ds}.

\begin{rem}[logarithmic monkey saddles of $\R V$]
The Hessian at the isolated parabolic points $\Log(\R V)$ vanishes.
Generic parabolic points sitting on hyperbolic components of $\Log(\R V)$
look like so-called monkey saddles (given in some local coordinates $(x,y,z)$
by $z=x(y^2-x^2)$).

Logarithmic monkey saddles do not appear on generic {\em smooth} surfaces in $(\R^*)^3$.
But they do appear on generic {\em real algebraic} surfaces in $(\R^*)^3$.
In particular, they appear on every torically maximal surface of sufficiently high degree.

The counterpart on the elliptic components of $\Log(\R V)$,
the {\em imaginary monkey saddles}, are locally given by $z=x(y^2+x^2)$.
\end{rem}
%}

\subsubsection{General case}
\label{hd}
Let $\R V\subset\rtor$ be a hypersurface.
Theorems \ref{kriv} and \ref{pov} have a weaker version
that holds for an arbitrary $n$.

\begin{thm}[\cite{M-r}]
If $\R V$ is torically maximal then every compact component of $\R V$ is a sphere.
All these $(n-1)$-spheres bound disjoint $n$-balls in $\rtor$.
\end{thm}

The following theorem is a counterpart of Theorem \ref{parpov} and a weaker version
of Theorem \ref{logkriv}.
\begin{thm}[\cite{M-r}]
The parabolic locus of $\Log(\R V)\subset\R^n$ is of codimension 2
if $\R V$ is torically maximal.
\end{thm}

Existence of torically maximal hypersurfaces for a given polyhedron
$\Delta$ seems to be a challenging question if $n>2$.

\ignore{
\subsection{Maximality conditions for non-Archimedean amoebas}
Let $\am_K\subset\R^n$ be a non-Archimedean amoeba (see \ref{nA})
whose Newton polyhedron is $\Delta$.
\begin{prop}
The number of vertices of $\am_K$ is not greater than $n!\Vol\Delta$.
\end{prop}
This proposition can be deduced from Theorem \ref{thmka} and the fact
that the smallest possible volume of a convex lattice polyhedron is $\frac{1}{n!}$.

\begin{defn}
\label{nAmax}
A non-Archimedean amoeba $\am_K$ is called {\em maximal} if the number of its
vertices equals to $n!\Vol\Delta$.
\end{defn}

\begin{rem}
\label{remnAmax}
For some choices of $\Delta$ maximal amoebas do not exist.
We can take, for instance, $\Delta\subset\R^3$ to be the tetrahedron
with vertices $(1,0,0)$, $(0,1,0)$, $(0,0,1)$ and $(0,0,k)$. Any valuation
on the corresponding coefficients would have to be linear. Its Legendre
transform would have just one vertex while $\Vol\Delta=\frac{k}{n!}$.
Note that these polyhedra were used by B. Bertrand to show
absence of real maximal surfaces, see Remark \ref{bertrand}.

Nevertheless, if the toric variety corresponding to $\Delta$ is a projective space
or a product of projective spaces then maximal non-Archimedean amoebas exist.
This statement is implicitly contained in \cite{IV}.
\end{rem}
}

\subsection{Amoebas and dimers}
Amoebas and, in particular, the amoebas of simple Harnack curves have appeared in
a recent work of Kenyon, Okounkov
and Sheffield on dimers, see \cite{KOS} and \cite{KO}.
In particular, Figure 1 of \cite{KOS} sketches
a probabilistic approximation of the amoeba of a line in the plane.

\begin{figure}[h]
\centerline{
\psfig{figure=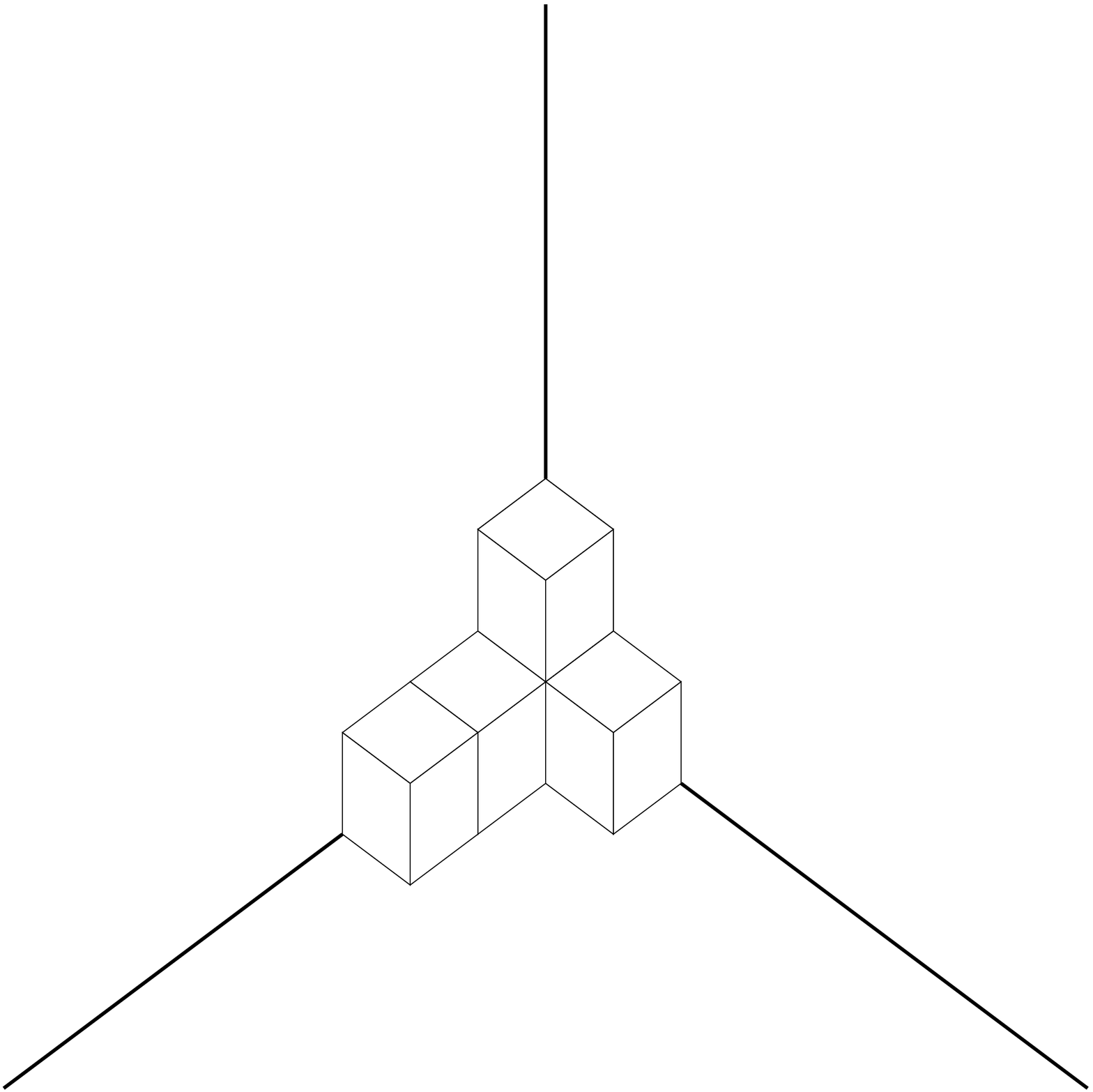,height=2in,width=2in}\hspace{.5in}
\psfig{figure=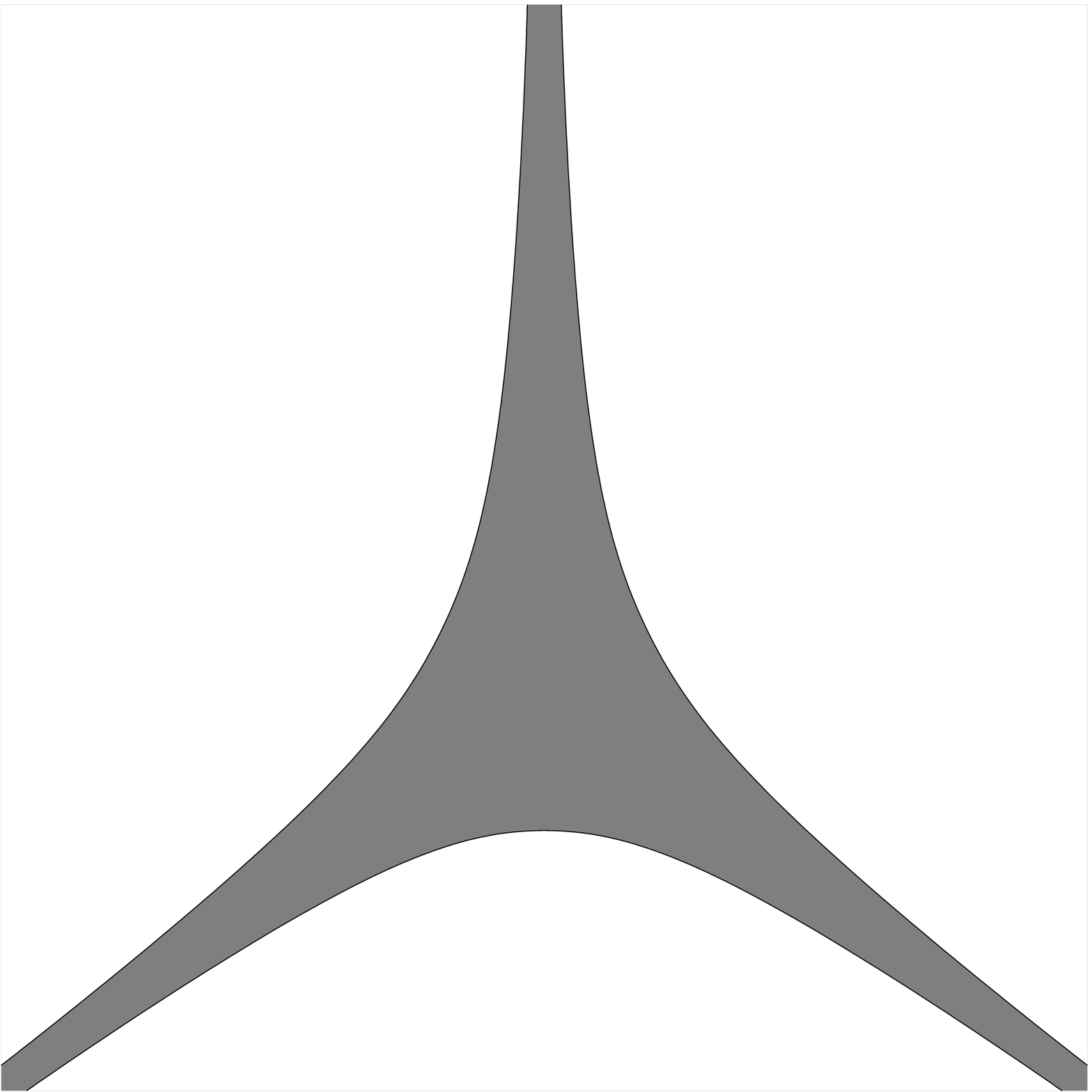,height=2in,width=2in}}
\caption{\label{fan} The fan with the dimer tiling
and the statistical accumulation $R$.}
\end{figure}

One starts from the negative octant
$$O=\{(x,y,z)\in\R^3\ |\ x<0, y<0, z<0\}.$$
Its projection onto $\R^2$ along the vector $(1,1,1)$
defines a fan with 3 corners, see Figure \ref{fan}.
For each $(x_0,y_0,z_0)\in\R^3$ let
$$Q_{(x_0,y_0,z_0)}=\{(x,y,z)\in\R^3\ |\
x_0-1<x\le x_0, y_0-1<y\le y_0, z_0-1<z\le z_0\}$$
be the unit cube with the ``outer" vertex $(x_0,y_0,z_0)$.
Let us fix a large natural number $N$ and remove $N$ such unit cubes from $O$
according to the following procedure.

At the first step we remove $Q_{(0,0,0)}$.
The region $O\setminus Q_{(0,0,0)}$ has three outer vertices,
namely $(-1,0,0)$, $(0,-1,0)$ and $(0,0,-1)$.
At the second step we remove a unit cube whose outer vertex is
one of these three and proceed inductively.
For each $N$ we have a finite number of possible resulting regions $O'$.
%for $O$ after removing the $N$ cubes.
The projection of such region defines a tiling by diamond-shaped figures
(dimers) as in Figure \ref{fan}. Clearly there is no more than $3N$
dimers in the tiling.
Each dimer in $\R^2$ is assigned a weight
in a double-periodic fashion with some integer period vectors.
The probability of a tiling is determined by these weights.

%To each such resulting region can be prescribed its own probability
%depending on the cubes that were removed.
%According to \cite{KOS} if the assigned probability satisfies
%to a certain periodicity condition then a certain scaling of $O'$
It is shown in \cite{KOS} that after some rescaling the union
of the dimer tiles converges to some limiting region $R\subset\R^2$
that depends only on the choice of the (periodic) choice
of weights of the dimers when $N\to\infty$.

Furthermore, according to \cite{KOS} there exists
a simple Harnack curve $V$ with the amoeba $\am\subset\R^2$
such that $$R=T(\am)$$
%from the overgraph $O\Gamma(N_f)$
%of the Ronkin function $N_f$ as
%$$R=T(O\Gamma(N_f)),$$
for the linear transformation $T
=\begin{pmatrix}
\frac{\sqrt{3}}{2} & -\frac{\sqrt{3}}{2}\\
\frac12 & \frac12
\end{pmatrix}$ in $\R^2$.
The curve $V$ is a line (as in Figure~\ref{fan})
if all the dimer weights are the same.
For other periodic weight choices any simple Harnack curve can appear.

Using such dimer interpretation Kenyon and Okounkov \cite{KO}
have constructed an explicit parameterization for the set
of all simple Harnack curves of the same degree.
It is shown in \cite{KO} that this set
is contractible.

\part{TROPICAL GEOMETRY}
\section{Tropical degeneration and the limits of amoebas}
\subsection{Tropical algebra}
\begin{defn}
The {\em tropical semifield} $\Rtr$ is the set of real numbers $\R$
equipped with the following two operation called {\em tropical addition}
and {\em tropical multiplication}. We use quotation marks to
distinguish tropical arithmetical operations from the standard ones.
For $x,y\in\Rtr$ we set $``x+y"=\max\{x,y\}$ and $``xy"=x+y$.
\end{defn}

This definition appeared in Computer Science.
The term ``tropical" was given in honor of
Imre Simon who resides in S\~ao Paolo, Brazil (see \cite{Pin}).
Strictly speaking, the tropical addition in Computer Science
is usually taken to be the minimum (instead of the maximum),
but, clearly, the minimum generates an isomorphic semifield.

The semifield $\Rtr$ lacks the subtraction. However it is not
needed to define polynomials.
Indeed the tropical polynomial is defined as
$$``\sum\limits_j a_jx^j"=\max\limits_j <j,x>+a_j$$
for any finite collections of coefficients $a_j\in\Rtr$
parameterized by indices $j=(j_1,\dots,j_n)\in\Z^n$.
Here $x=(x_1,\dots,x_n)\in\R^n$, $x^j=x_1^{j_1}\dots x_n^{j_n}$
and $<j,x>=j_1x_1+\dots+j_nx_n$.

Thus the tropical polynomials are piecewise-linear functions.
They are simply the Legendre transforms of the function
$j\mapsto -a_j$ (this function is defined only on finitely many points,
but its Legendre transform is defined everywhere on $\R^n$).

It turns out that these polynomials are responsible for some
piecewise-linear geometry in $\R^n$ that is similar in many ways
to the classical algebraic geometry defined by the polynomials
with complex coefficients. Furthermore, this tropical geometry
ca be obtained as the result of a certain degeneration of the
(conventional) complex geometry in the torus $\tor$.

\subsection{Patchworking as tropical degeneration}
In 1979 Viro discovered a {\em patchworking} technique for
construction of real algebraic hypersurfaces, see \cite{V}.
Fix a convex lattice polyhedron $\Delta\in\R^n$.
Choose a function $v:\Delta\cap\Z^n\to\R$.
The graph of $v$ is a discrete set of points in $\R^n\times\R$.
The overgraph is a family of parallel rays. Thus the convex hull
of the overgraph is a semi-infinite polyhedron $\tilde\Delta$.
The facets of $\tilde\Delta$ which project isomorphically to
$\R^n$ define a subdivision of $\Delta$ into smaller
convex lattice polyhedra $\Delta_k$.

Let $F(z)=\sum\limits_{j\in\Delta}a_jz^j$ be a generic polynomial in the class of
polynomial whose Newton polyhedron is $\Delta$.
The {\em truncation} of $F$ to $\Delta_k$ is $F_{\Delta_k}=\sum\limits_{j\in\Delta_k}a_jz^j$.
The {\em patchworking polynomial} $f$ is defined by formula
\begin{equation}
\label{patch}
f^v_t(z)=\sum\limits_j a_jt^{v(j)}z^j,
\end{equation}
$z\in\R^n$, $t>1$ and $j\in\Z^n$.

Consider the hypersurfaces $V_{\Delta_k}$ and $V_t$ in $\tor$ defined
by $F_{\Delta_k}$ and $f^v_t$.
If $F$ has real coefficients then we denote $\R V_{\Delta_k}=V_{\Delta_k}\cap\rtor$
and $\R V_t=V_t\cap\rtor$.
Viro's patchworking theorem \cite{V} asserts that for large values of $t$
the hypersurface $\R V_t$ can be obtained from $\R V_{\Delta_k}$ by a certain patchworking
procedure. The same holds for amoebas of the hypersurfaces $V_t$ and $\R V_{\Delta_k}$.
In fact patchworking of real hypersurfaces can be interpreted as the real version
of patchworking of amoebas (cf. Appendix in \cite{Mi}).
%Below we describe a special case of amoeba's patchworking
%in terms of the so-called dequantization.
%
%\subsection{Maslov's dequantization}
%\label{deq}
It was noted by Viro in \cite{V2000} that patchworking is related to so-called
{\em Maslov's dequantization} of positive real numbers.

Recall that a {\em quantization} of a semiring $R$ is a family of semirings $R_h$,
$h\ge 0$ such that $R_0=R$ and $R_t\approx R_s$ as long as $s,t>0$, but
$R_0$ is not isomorphic to $R_t$.
The semiring $R_h$ with $h>0$ is called a {\em quantized} version of $R_0$.

Maslov (see \cite{Ma}) observed that the ``classical" semiring $\R_+$
of real positive number is a quantized version of some other ring in this sense
Let $R_h$ be the set of positive numbers with the usual multiplication
and with the addition operation
$z\oplus_h w=(z^{\frac{1}{h}}+w^{\frac{1}{h}})^h$
for $h>0$ and
$z\oplus_h w=\max\{z,w\}$
for $h=0$.
Note that
$$\lim\limits_{h\to 0} (z^{\frac{1}{h}}+w^{\frac{1}{h}})^h=\max\{z,w\}$$
and thus this is a continuous family of arithmetic operations.

The semiring $R_1$ coincides with the standard semiring $\R_+$.
The isomorphism between $\R_+$ and $R_h$ with $h>0$ is given
by $z\mapsto z^h$. On the other hand the semiring $R_0$ is not
isomorphic to $\R_+$ since it is idempotent, indeed $z+z=\max\{z,z\}=z$.

%\subsection{Logarithmic dequantization}
%\label{logdeq}
Alternatively we may define the dequantization deformation with the help
of the logarithm.
The logarithm $\log_t$, $t>1$, induces a semiring structure on $\R$ from $\R_+$,
$$x\oplus_t y=\log_t(t^x+t^y),\ x\otimes_t y=x+y,\ x,y\in\R.$$
Similarly we have $x\oplus_\infty y=\max\{x,y\}$.
Let $R^{\log}_t$ be the resulting semiring.
\begin{prop}
The map $\log:R_h\to R^{\log}_t$, where $t=e^{\frac1h}$, is an isomorphism.
\end{prop}

%\subsection{Patchworking as a dequantization}
The patchworking polynomial \eqref{patch} can be viewed as
a deformation of the polynomial $f^v_1$. We define a similar
deformation with the help of Maslov's dequantization.
Instead of deforming the coefficients we keep the coefficients
but deform the arithmetic operations.

Choose any coefficients $\alpha_j$, $j\in\Delta$.
Let $\phi_t:(R^{\log}_t)^n\to R^{\log}_t$, $t\ge e$, be a polynomial whose
coefficients are $\alpha$, i.e. $$\phi_t(x)=\bigoplus_t(\alpha_j+jx),\ x\in\R^n.$$
Let $\Log_t:\tor\to\R^n$ be defined by $(x_1,\dots,x_n)=(\log|z_1|,\dots,\log|z_n|)$.
\begin{prop}[Maslov \cite{Ma},Viro \cite{V2000}]
The function $f_t=(\log_t)^{-1}\circ\phi_t\circ\Log_t:(\R_+)^n\to\R_+$ is a polynomial with respect to the standard
arithmetic operations in $\R_+$, namely we have
$$f_t(z)=\sum\limits_j t^{\alpha_j} z^j.$$
\end{prop}

This is a special case of the patchworking polynomial \eqref{patch}. The coefficients
$\alpha_j$ define the function $v:\Delta\cap\Z^n\to\R$.

\subsection{Limit set of amoebas}
Let $V_t\subset\tor$ be the zero set of $f_t$ and let $\am_t=\Log_t(V_t)\subset\R^n$.
Note that $\am_t$ is the amoeba of $V_t$ scaled $\log t$ times.
Note also that the family $f_t=\sum\limits_j t^{\alpha_j} z^j$ can be considered
as a single polynomial whose coefficients are powers of $t$.
Such coefficients are a very simple instance
of the so-called {\em Puiseux series}.

The field $K$ of the real-power Puiseux series is obtained
from the field of the Laurent series in $t$ by taking
the algebraic closure first and then taking the metric completion
with respect to the ultranorm
$$||\sum a_j t^j||=\min\{ j\in\R \ |\ a_j\neq 0\}.$$
The logarithm $\val:K^*\to\R$ of this norm is
an example of the so-called non-Archimedean valuation
as $\val(a+b)\le\max\{\val(a)+\val(b)\}$ and
$\val(ab)=\val(a)+\val(b)$ for any $a,b\in K^*=K\setminus\{0\}$.

\begin{defn}[Kapranov \cite{Ka}]
Let $V_K\subset (K^*)^n$ be an algebraic variety.
Its {\em (non-Archimedean) amoeba} is
$$\am_K=\Val(V_k)\subset\R^n,$$
where $\Val(z_1,\dots,z_n)=(\val(z_1),\dots,\val(z_n))$.
\end{defn}

%\begin{thm}{\cite{Ka}}
%If $V_K\subset (K^*)^n$ is a hypersurface given
%as the zero set of a polynomial $f:K^n\to K$
%then
%\end{thm}

We have a uniform convergence of the addition operation
in $R^{\log}_t$ to the addition operation in $R^{\log}_\infty$.
As it was observed by Viro it follows from the following inequality
$$\max\{x,y\}\le x\oplus_t y=\log_t(t^x+t^y)\le \max\{x,y\}+\log_t 2.$$
More generally, we have the following lemma.
\begin{lem}\label{logestimate}
$$\max\limits_{j\in\Delta}(\alpha_j+jx) \le \phi_t(x)\le \max\limits_{j\in\Delta}(\alpha_j+jx)+
\log N,$$
where $N$ is the number of lattice points in $\Delta$.
\end{lem}

Recall that the {\em Hausdorff metric} is defined on closed subsets
$A,B\subset\R^n$ by
$$d_{\operatorname{Hausdorff}}(A,B)=\max\{\sup_{a\in A} d(a,B),
\sup_{b\in B} d(A,b)\},$$
where $d$ is the Euclidean distance in $\R^n$.
The following theorem is a corollary of Lemma \ref{logestimate}.
\begin{thm}[Mikhalkin \cite{M-pp}, Rullg{\aa}rd \cite{R2}]
\label{amlim}
The subsets $\am_t\subset\R^n$ tend in the Hausdorff metric to $\am_K$
when $t\to 0$.
\end{thm}
Recall that in our setup $t>0$. Alternatively
we may replace $t$ with $\frac{1}{t}$ to
get a limit with $t\to +\infty$.

%Note that by Theorem \ref{thmka} $\am_K$ is obtained by patchworking of
%he amoebas of the truncations of $f_e$ to smaller polyhedra $\Delta_k$ (see \ref{pp}).

\subsection{Tropical varieties and non-Archimedean amoebas}
We start by defining tropical hypersurfaces.
The semiring $\Rtr$ lacks (additive) zero so the tropical hypersurfaces
are defined as singular loci and not as zero loci.
Let $F:\R^n\to\R$ be a tropical polynomial.
It is a continuous convex piecewise-linear function.
Unless $F$ is linear it is not everywhere smooth.
\begin{defn}
{\em The tropical variety $V_F\subset\R^n$ of $F$} is the set of all
points in $\R^n$ where $F$ is not smooth.
\end{defn}
Equivalently we may define $V_F$ as the set of points where
more than one monomial of $F(x)=``\sum a_jx^j"$ reaches the maximum.

Let us go back to the non-Archimedean field $K$ of Puiseux series.
Let $$f(z)=\sum\limits_j \alpha_jz^j,$$ $\alpha_j\in K$,
$j\in\Z^n$, $z\in K^n$, be a polynomial that defines a hypersurface
$V_K\subset (K^*)^n$ and let $\am_K\subset\R^n$ be the corresponding
non-Archimedean amoeba. We form a tropical polynomial
$$F(x)=\sum\limits_j \val(\alpha_j)x^j,$$
$x\in\R^n$.

Kapranov's description \cite{Ka} of the non-Archimedean amoebas
can be restated in the following way.
\begin{thma}[\cite{Ka}]
The amoeba $\am_K$ coincides with the tropical hypersurface $V_F$.
\end{thma}

Definition of tropical varieties in higher codimension in $\R^n$
gets somewhat tricky as intersections of tropical hypersurfaces are
not always tropical.
As is was suggested in \cite{RGST}
non-Archimedean amoebas provide a byway for such definition as tropical
varieties can be simply defined as non-Archimedean amoebas
for algebraic varieties in $(K^*)^n$.

In the next section we concentrate on the study of tropical curves.
References to some higher-dimensional tropical varieties treatments
include \cite{M-pp} for the case of hypersurfaces and
\cite{SS} for the case of the Grassmanian varieties.

\section{Calculus of tropical curves in $\R^n$}
\subsection{Definitions}
Let $\bar\Gamma$ be a finite graph whose edges are weighted by
natural numbers.
Let ${\mathcal V}_1$ be the set of 1-valent vertices of $\Gamma$.
We set
$$\Gamma=\bar\Gamma\setminus {\mathcal V_1}.$$

\begin{defn}[Mikhalkin \cite{M-en}]
%, cf. Aharony-Hanany-Kol \cite{AHK}]
\label{tropcur} A proper map $h:\Gamma\to\R^n$ is called {\em a
parameterized tropical curve} if it satisfies to the following two
conditions.
\begin{itemize}
\item For every edge $E\subset\Gamma$ the restriction $h|_E$ is an
embedding. The image $h(E)$ is contained in a line $l\subset\R^n$
such that the slope of $l$ is rational.
\item For every vertex
$V\in\Gamma$ we have the following property. Let
$E_1,\dots,E_m\subset\Gamma$ be the edges adjacent to $V$, let
$w_1,\dots,w_m\in\N$ be their weights and let
$v_1,\dots,v_m\in\Z^n$ be the primitive integer vectors from $V$
in the direction of the edges. We have
\begin{equation}
\label{balance}
\sum\limits_{j=1}^m w_jv_j=0.
\end{equation}
\end{itemize}
Two parameterized tropical curves $h:\Gamma\to\R^n$
and $h':\Gamma'\to\R^n$ are called {\em equivalent} if there exists a
homeomorphism $\Phi:\Gamma\to\Gamma'$ which respects the weights
of the edges and such that $h=h'\circ\Phi$.
We do not distinguish equivalent parameterized tropical curves.

The image $$C=h(\Gamma)\subset\R^n$$ is called the (unparameterized)
tropical curve. It is a weighted piecewise-linear graph in $\R^n$.
Note that the same curve $C\subset\R^2$ may admit non-equivalent
parameterizations.
The curve $C$ is called {\em irreducible} if $\Gamma$ is connected
for any parameterization. Otherwise the curve is called reducible.
\end{defn}

\begin{figure}[h]
\centerline{
\psfig{figure=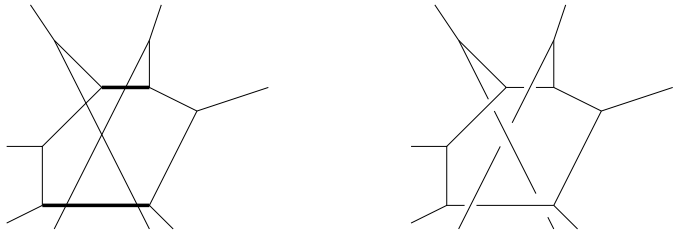,height=1.1875in,width=3.5in}}
\caption{\label{trcr23} A tropical curve in $\R^2$ and its possible lift
to $\R^3$. The edges of weight 2 are bold (at the left picture).
Note that lifts of such edges
can have weight 1.}
\end{figure}

\begin{rmk}
In dimension 2 the notion of tropical curve coincides with
the notion of $(p,q)$-webs introduced by Aharony, Hanany and Kol
in \cite{AHK} (see also \cite{AH}).
\end{rmk}

It is convenient to prescribe a multiplicity to a 3-valent
vertex $A\in \Gamma$ of the tropical curve $h:\Gamma\to\R^n$
as in \cite{M-en}.
As in Definition \ref{tropcur} let
$w_1,w_2,w_3$ be their weights of the edges of $h(\Gamma)$
adjacent to $A$ and let
$v_1,v_2,v_3$ be the primitive integer vectors in the direction of the edges.
\begin{defn}\label{multvert}
The {\em multiplicity} of a 3-valent vertex $A$ in $h(\Gamma)$
is $w_1w_2|v_1\times v_2|$.
Here $|v_1\times v_2|$ is the ``length of the vector product of $v_1$
and $v_2$" in $\R^n$ being interpreted as the area of the parallelogram spanned
by $v_1$ and $v_2$.
Note that
$$w_1w_2|v_1\times v_2|=w_2w_3|v_2\times v_3|=w_3w_1|v_3\times v_1|$$
since $v_1w_1+v_2w_2+v_3w_3=0$ by Definition \ref{tropcur}.
%The (global) {\em multiplicity of a simple curve} is the product of multiplicities
%of all its vertices.
\end{defn}

\begin{figure}[h]
\centerline{\psfig{figure=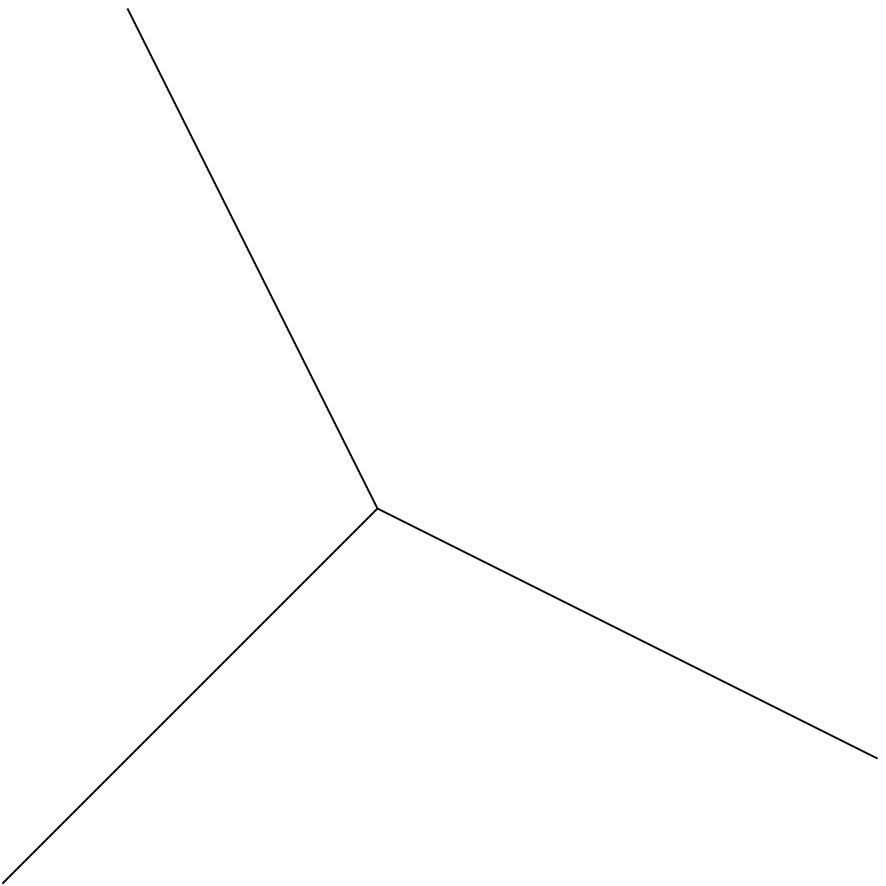,height=1in,width=1in}\hspace{1in}
\psfig{figure=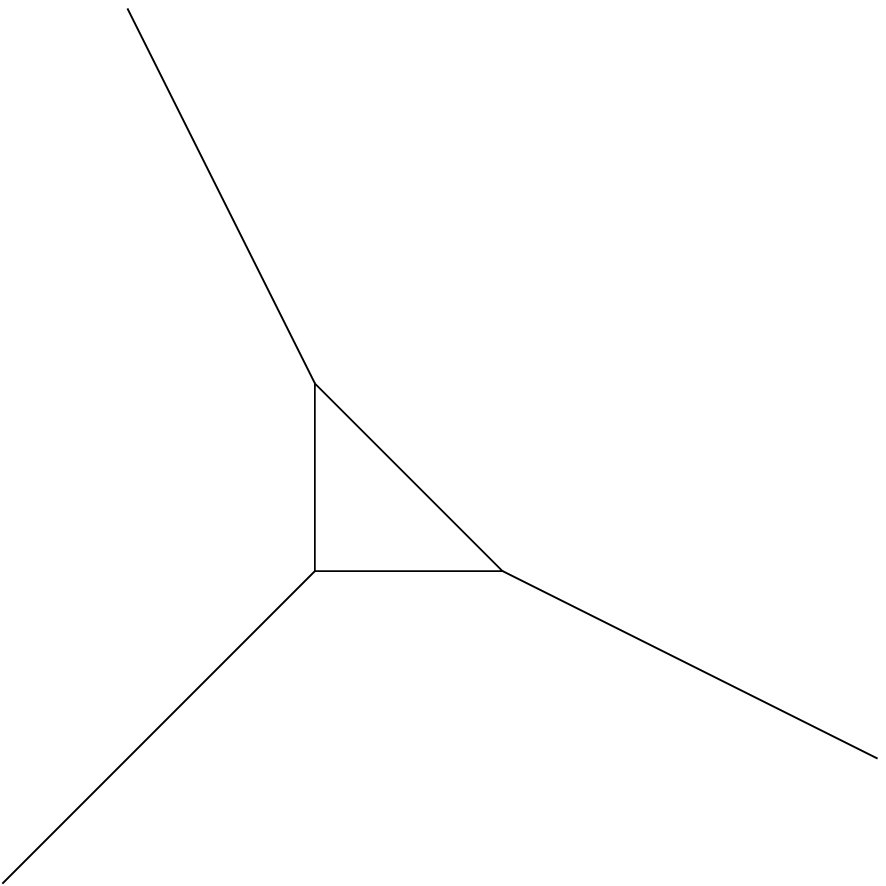,height=1in,width=1in}}
\caption{\label{3m3} Deformation of a multiple 3-valent vertex}
\end{figure}

If the multiplicity of a vertex is greater than 1 then it
is possible to deform it with an appearance of a new cycle
as in Figure \ref{3m3}.

\subsection{Degree, genus and the tropical Riemann-Roch formula}
Heuristically, the degree of a tropical curve $C\subset\R^n$ is the set of
its asymptotic directions. For each end of a tropical curve $C=h(\Gamma)$
we fix a primitive integer vector parallel to this ray in the outward direction
and multiply it by the weight of the corresponding (half-infinite) edge.
Doing this for every end of $C$ we get a collection ${\mathcal C}$
of integer vectors in $\Z^n$.

Let us add all vectors in ${\mathcal C}$
that are positive multiples of each other.
The result is a set
${\mathcal T}=\{\tau_1,\dots,\tau_q\}\subset\Z^n$
of non-zero integer vectors such that $\sum\limits_{j=1}^q\tau_j=0$.
Note that in this set we do not have positive multiples
of each other, i.e. if $\tau_j=m\tau_k$ for $m\in\N$ then $\tau_j=\tau_k$.

\begin{defn}[\cite{M-en}]\label{trdeg}
The set ${\mathcal T}$ is called {\em the degree of the tropical curve $C\subset\R^n$}.
{\em The genus} of a parameterized tropical curve $h:\Gamma\to\R^n$
is $\dim(H_1(\Gamma))+1-\dim(H_0(\Gamma))$ so that if $\Gamma$ is connected
then it coincides with the number of cycles $\dim(H_1(\Gamma))$ in $\Gamma$.
The genus of a tropical curve $C\subset\R^n$ is the minimal genus
among all the parameterization $C=h(\Gamma)$.
\end{defn}

There is an important class of tropical curves that behaves
especially nice with respect to a genus-preserving deformation.
\begin{defn}[\cite{M-en}]\label{simplecurve}
A parameterized tropical curve $h:\Gamma\to\R^n$ is called {\em simple}
if
\begin{itemize}
\item $\Gamma$ is 3-valent,
\item $h$ is an immersion,
\item if $a,b\in\Gamma$ are such that $h(a)=h(b)$ then neither $a$
nor $b$ can be a vertex of $\Gamma$.
\end{itemize}
In this case the image $h(\Gamma)$ is called
a {\em a simple tropical curve}.
%(or {\em immersed 3-valent})
\end{defn}

Simple curves locally deform in a linear space.
%It turns out that as in the classical case
%there
\begin{thm}[Tropical Riemann-Roch, \cite{M-en}]\label{trRR}
Let $h:\Gamma\to\R^n$
be a simple tropical curve, where $\Gamma$ is a graph with $x$ ends.
Non-equivalent tropical curves of the same genus and with the
same number of ends close to $h$ locally
form a $k$-dimensional real vector space, where
$$k\ge x+(n-3)(1-g).$$
\end{thm}
If the curve is non-simple then its space of deformation is locally
piecewise-linear.

\subsection{Enumerative tropical geometry in $\R^2$}\label{secdim2}
We start by considering the so-called
"curve counting problem" for the complex torus $\tordva$.

Any algebraic curve $V\subset\tordva$ is defined by a polynomial
$$f(z,w)=\sum\limits_{j,k} a_{jk}z^jw^k.$$
Recall that from the topological viewpoint
the degree of a variety is its homology class
in the ambient variety.
Here we have a difficulty caused by non-compactness
of $\tordva$.

Help is provided by {\em the Newton polygon}
$$\Delta(f)=\operatorname{Convex Hull}\{
(j,k)\ |\ a_{jk}\neq 0 \}$$
of $f$.
The polygon $\Delta=\Delta(f)$ can be interpreted as {\em the (toric) degree}
of $V$. Indeed being a compact lattice polygon
$\Delta$ defines a compact toric surface $\C T_{\Delta}\supset\tordva$,
e.g. by taking the closure of the image under {\em the Veronese embedding}
$\tordva\to\cp^{\#(\Delta\cap\Z^2)}$ (see e.g. \cite{GKZ}).
The closure of $V$ in $\C T_\Delta$ defines a homology class
induced from the hyperplane section by the Veronese embedding.

Note that the definition of the toric degree agrees with
its tropical counterpart in Definition \ref{trdeg}.
Indeed, for each side $\Delta'$ of $\Delta$ we can take
the primitive integer normal vector in the outward direction and multiply
it by the lattice length $\#(\Delta'\cap\Z^2)-1$ of the side.
The result is a tropical degree set ${\mathcal T}(\Delta)$.
Accordingly we define
$$x=\#(\dd\Delta\cap\Z^2)$$
which is the number of ends of a general curve of degree $\Delta$
in $\tordva$.

An irreducible curve $V$ has {\em geometric genus} which
is the genus of its normalization $\tilde{V}\to V$.
In the case when $V$ is not necessarily irreducible
it is convenient to define the genus as the sum
of the genera of all irreducible components minus the number of
such components plus one.

Let us fix the genus (i.e. a number $g\in\Z$)
and the toric degree (i.e. a polygon $\Delta\subset\R^2$).
Let $$\ppp=\{p_1,\dots,p_{x+g-1}\}\subset\tordva$$
be an configuration of $x+g-1$ general points in $\tordva$.
We set $\nclr(g,\Delta)$ to be equal to be the number of
curves in $\tordva$ of genus $g$ and degree $\Delta$
passing through $\ppp$. Similarly we set $\ncl(g,\Delta)$
to be the number of irreducible curves among them.

These numbers are close relatives of the Gromov-Witten
invariants of $\C T_\Delta$ (see \cite{KM} for the definition).
In the case when $\C T_\Delta$ is smooth Fano they
coincide with the corresponding Gromov-Witten invariants.
The numbers $\nclr(g,\Delta)$ and $\ncl(g,\Delta)$ have
tropical counterparts.

For a fixed genus $g$ and a toric degree $\Delta$ we fix
a configuration $$\qqq=\{r_1,\dots,r_{x+g-1}\}\subset\R^2$$
of $x+g-1$ general points in the tropical plane $\R^2$
(for a rigorous definition of tropical general position see \cite{M-en}).
We have a finite number of tropical curves of genus $g$ and degree
$\Tau(\Delta)$ passing through $\qqq$, see \cite{M-en}.
Generically all such curves are simple (see Definition \ref{simplecurve}.
However unlike the situation in $\tordva$
the number of such curves is different for different configurations
of $x+g-1$ general point.

\begin{defn}[\cite{M-en}]\label{multdim2}
The multiplicity $\mult(h)$ of a simple
tropical curve $h:\Gamma\to\R^2$ of degree $\Delta$ and
genus $g$ passing via $\qqq$ equals to the product of
the multiplicities of the (3-valent) vertices of $\Gamma$.
(see Definition \ref{multvert}).
\end{defn}

\begin{thm}[\cite{M-en}]\label{thm2dim}
The number of irreducible tropical curves of genus $g$ and degree $\Delta$
passing via $\qqq$ and counted with multiplicity from Definition \ref{multdim2}
equals to $\ncl(g,\Delta)$.

The number of all tropical curves of genus $g$ and degree $\Delta$
passing via $\qqq$ and counted with multiplicity from Definition \ref{multdim2}
equals to $\nclr(g,\Delta)$.
\end{thm}

\begin{exa}
Figure \ref{picdeg3} shows a (generic) configuration of 8 points
$\qqq\subset\R^2$ and all curves of genus 0 and of projective degree 3 passing
through $\qqq$. Out of these nine curves eight have multiplicity 1 and one
(with a weight 2 edge) has multiplicity 4. All the curves are irreducible.
Thus $\ncl(g,\Delta)=\nclr(g,\Delta)=12$.
\begin{figure}[h]
\centerline{\psfig{figure=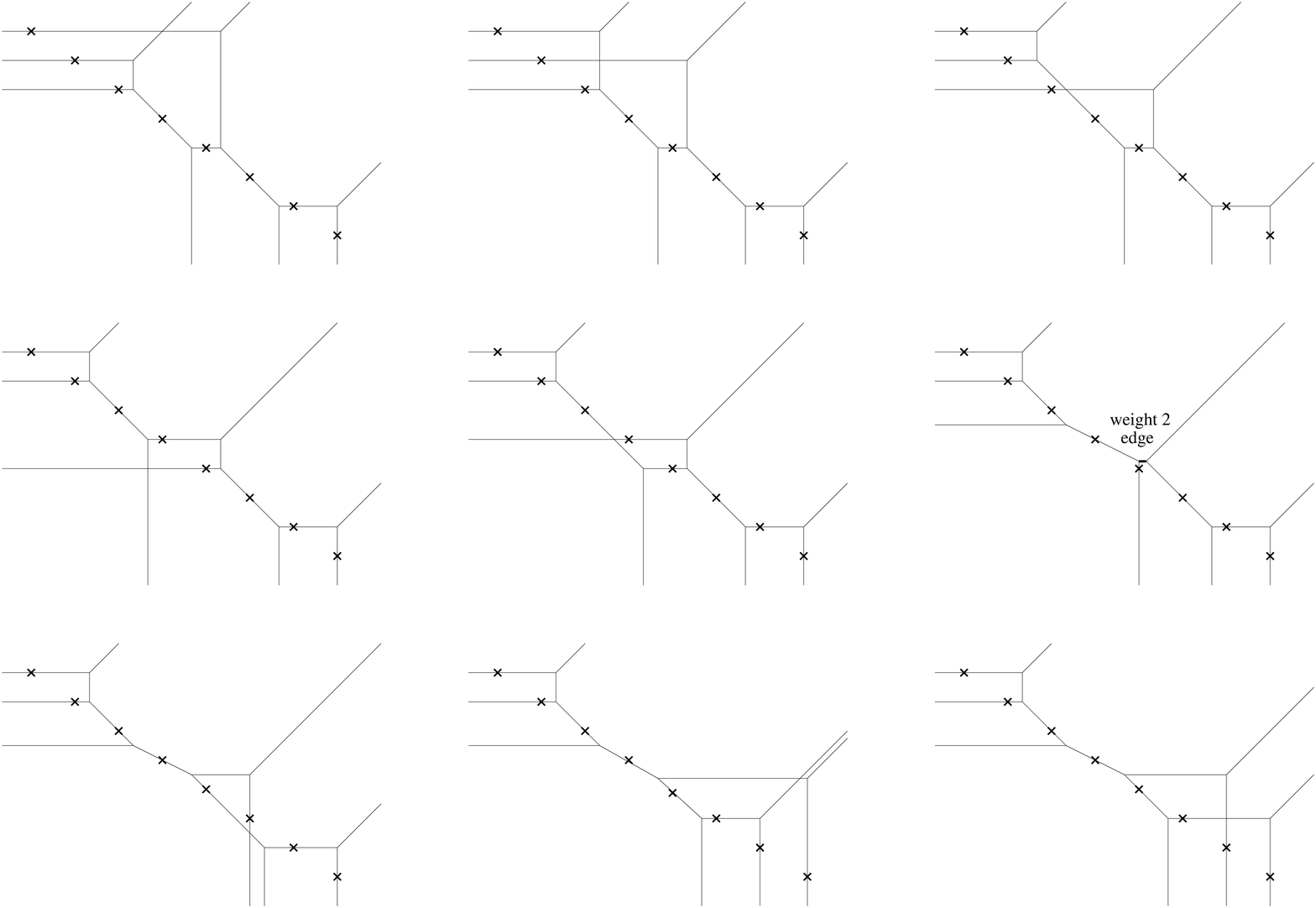,height=3.3in,width=4.84in}}
\caption{\label{picdeg3} Tropical projective rational cubics via 8 points.}
\end{figure}
\end{exa}

Theorem \ref{thm2dim} thus reduces the problem of finding
$\ncl(g,\Delta)$ and $\nclr(g,\Delta)$ to the corresponding
tropical problems. Furthermore, it allows to use any general
configuration $\qqq$ in the tropical plane $\R^2$
(as it implies that the answer is independent of $\qqq$).
We can take the configuration $\qqq$ on the same affine
(not tropical) line $L\subset\R^2$ and still insure
tropical general position as long as the slope of $L$
is irrational. It was shown in \cite{M-en}
that such curves are encoded by lattice paths of
length $x+g-1$ connecting a pair of vertices in $\Delta$.

Namely, the slope of $L$ determines a linear function
$\lambda:\R^2\to\R$ such that $\lambda|_{\Delta\cap\Z^2}$
is injective and thus a linear order on the lattice points
of $\Delta$. There is a combinatorial rule (see \cite{M-cr}
or \cite{M-en}) that associates a non-negative integer
multiplicity to every $\lambda$-increasing lattice path
of length $x+g-1$, i.e. to every order-increasing sequence
of lattice points of $\Delta$ that contains $x+g$ points.
This multiplicity is only non-zero if the first and the
last points of the sequence are the points where $\lambda|_{\Delta}$
reaches its minimum and maximum.

\begin{exa}
The tropical curves from Figure \ref{picdeg3} are described
by the lattice paths from Figure \ref{pathdeg3} shown
together with their multiplicities.
Here the first path describes the first 3 tropical curve
from Figure \ref{picdeg3}, the second --- the next two paths,
the third --- the next curve (which itself corresponds to 4
distinct holomorphic curves), the fourth --- the next curve and
the fifth --- the last two tropical curves from Figure \ref{picdeg3}.
These paths are $\lambda$-increasing for $\lambda(x,y)=y-(1+\epsilon x)$,
where $\epsilon>0$ is very small.
\begin{figure}[h]
\centerline{\psfig{figure=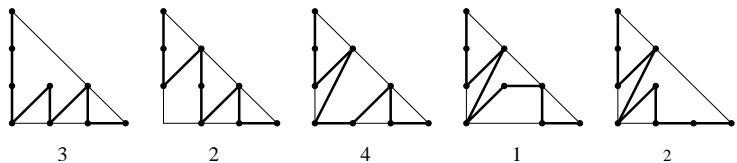,height=.8in,width=3.8in}}
\caption{\label{pathdeg3} The lattice paths describing the tropical
curves from Figure \ref{picdeg3} and the path multiplicities.}
\end{figure}
\end{exa}

\subsection{Enumerative tropical geometry in $\R^3$ (and higher dimension)}
The results of previous subsections can be established with the help of
the following restatement of Theorem \ref{amlim} in the case of $\R^2$.
\begin{lem}
If $C=h(\Gamma)\subset\R^2$ is a tropical curve then there exists
a family $V_t\subset\tordva$ of holomorphic curves for $t>0$ such
that $\Log_t(V_t)=C$. Here the degree of $C$ coincides with the
degree of $V_t$.
\end{lem}

The situation is more complicated if $n>2$
as such statement is no longer true for {\em all}
tropical curves in $\R^n$.
\begin{exa}
Consider the graph $C'\subset\R^2\subset\R^3$ depicted on
%(part of) the planar projective cubic depicted on
Figure \ref{supercubic}. This set can be obtained by removing
three rays from a planar projective cubic curve.
Let $q_1,q_2,q_3\in \R^2$ be the end points of these rays.
Consider the curve
$$C=C'\cup\bigcup\limits_{j=1}^n
(\{(q_j,t)\ | t\le 0\}\cup\{(q_j+t,t)\ | t\ge 0\}).$$
It is easy to check that $C\subset\R^3$ is
a (spatial) projective curve of degree 3 and genus 1.
Suppose that $q_1,q_2,q_3$ are not tropically collinear,
i.e. are not lying on the same tropical line in $\R^2$
(e.g. we may choose $q_1,q_2,q_3$ to be in tropically
general position). Then $C$ cannot be obtained as the limit
of $\Log_t(V_t)$ for cubic curves $V_t\subset\cp^3$
(since $\Log_t(V_t)$ is not everywhere defined
$\Log_t(V_t)$ stands for $\Log_t(V_t\cap\tortri$).
\begin{figure}[h]
\centerline{\psfig{figure=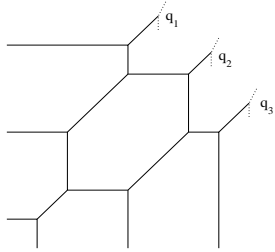,height=1.3in,width=1.4in}}
\caption{\label{supercubic} A planar part of a superabundant
spatial cubic}
\end{figure}

Indeed, any cubic curve $V_t\subset\cp^3$ of genus 1
is planar, i.e. is contained in a plane $H_t\subset\cp^3$.
%The plane $H_t$ is unique (unless $V_t$ is reduced).
%By the compactness argument for tropical
It is easy to see (after passing to a subsequence
cf. \cite{M-pp} and \cite{M-en}) that
there has to exist a limiting set $H$ for $\Log_t(H_t)$.
%Clearly $C\subset H$.
Furthermore, $H$ is a tropical hypersurface in $\R^3$ whose
Newton polyhedron is contained in the polyhedron of a hyperplane.
Since $C\subset H$ we can deduce that $H$ has to be a hyperplane.
But then the intersection of $H$ with $\R^2$ is
(up to a translation in $\R^2$) a union
of the negative quadrant $\{(x,y)\ |\ x\le 0, y\le 0\}$
and the ray $\{(t,t)\ |\ t\ge 0$. The points $p_j$ have to sit
on the boundary of the quadrant (which is impossible unless
they are tropically collinear).
\end{exa}

Note that the tropical Riemann-Roch formula (Theorem \ref{trRR})
is a strict inequality for the curve $C$.
In accordance with the classical terminology such curves
are called {\em superabundant}. Conversely, a tropical curve is
called {\em regular} if the Riemann-Roch formula turns into equality.
It is easy to see that all rational curves are regular and that
the superabundancy of $C$ is caused by the cycle contained in
an affine plane in $\R^3$. Conjecturally all regular curves
are limits of the corresponding complex amoebas. Hopefully
the technique developed in the Symplectic Field Theory,
see \cite{EGH} and \cite{Bo} can help to verify this conjecture.

Let us formulate a tropical enumerative tropical problem in $\R^n$.
We fix the genus $g$ and the degree $\Tau=\{\tau_1,\dots,\tau_q\}\subset\Z^n$.
In addition we fix a configuration $\qqq$ which consists of
some points and some higher dimensional tropical varieties in $\R^n$
in general position. Let $k$ be the sum of the codimensions of
all varieties in $\qqq$. For each $\tau_j$ let $x_j\in\N$ be the
maximal integer that divides it. Let $x=\sum\limits_{j=1}^gx_j$.

If $k=x+(1-g)(n-3)$ then the expected number of tropical curves
of genus $g$ and degree $\Tau$ passing through $\qqq$ is finite.
However there may exist positive-dimensional families of superabundant
curves of genus $g$ and degree $\Tau$ through $\qqq$.

One way to avoid this (higher-dimensional) difficulty is
to restrict ourselves to the genus zero case.
In this case one can assign multiplicities to tropical rational
curves passing through $\qqq$ so that the total number of tropical
curves counted with these multiplicities agrees with the number of curves
in the corresponding complex enumerative problem
%In particular, for $n=3$ multiplicities concentrate
%Unlike the planar case where the multiplicities are the products
%of vertices multiplicities here we
(details are subject to a future paper).
%(see a forthcoming paper \cite{M-rn} for details).
\ignore{
Below we give the details for
the $n=3$ case. In this case $\qqq$ consists of a collection
of points and a collection of 1-cycles.

Let $\Gamma\subset\
If $\qqq$ is generic then it intersect
}

\subsection{Complex and real tropical curves}
Tropical curve $C\subset\R^n$ can be presented as
images $C=\Log(B)$ for certain objects $B\subset\tor$
called {\em complex tropical curves}. (Recall that
$\Log:\tor\to\R^n$ is the coordinatewise logarithm of
the absolute value.)
Let $H_t:\tor\to\tor$ be the self-diffeomorphism defined
by $H_t(z_1,\dots,z_n)=(|z_1|^{\log(t)-1}z^1,\dots,|z_n|^{\log(t)-1}z_n)$.
We have $\Log_t(z)=\Log(H_t(z))$.
\begin{defn}\label{defB}
The set $B\subset\tor$ is called a complex tropical curve
if it satisfies to the following condition.
\begin{itemize}
\item For every $x\in\R^n$ there exist
a neighborhood $U\ni x$
and a family $V_t\subset\tor$, $t>1$ of holomorphic curves
such that
$$B\cap\Log^{-1}(U)=\lim\limits_{t\to+\infty}(H_t^{-1}(V_t)\cap U),$$
where the limit is taken with respect to the Hausdorff metric.
\item For every open set $U\subset\R^n$ for every component $B'$
of $B\cap\Log^{-1}(U)$ there exists a tropical curve $C'\subset\R^n$
such that projection $\Log(B')=C'\cap U$.
\end{itemize}
\end{defn}
It is easy to see that for every open edge of $E\subset C$ the
inverse image $\Log^{-1}(E)\cap B$ is a disjoint union
of holomorphic cylinders.
%Thanks to the second condition of Definition \ref{defB} we
We can prescribe the weights to this cylinder so that
the sum is equal to the weight of $E$.
(In fact the second condition in Definition \ref{defB}
is needed only to insure that the cylinder weights
in different neighborhoods are consistent.)

Complex tropical curves can be viewed as curves ``holomorphic"
with respect to a (maximally) degenerate complex structure in $\tor$.
Consider a family of almost complex structures $J_t$
induced from the standard structure on $\tor$ by
the self-diffeomorphism $H_t$, $t>1$.
For every finite $t$ it is an integrable complex structure
(isomorphic to the standard one by $H_t$).
The curves $H^1(V_t)$ are $J_t$-holomorphic as long as $V_t$ is
holomorphic (with respect to the standard, i.e. $J_e$-holomorphic structure).
The limiting $J_\infty$-structure is no longer complex or almost complex,
but it is convenient to view $B$ as a
``$J_\infty$-holomorphic curve".

If $C=\Log(B)$ admits a parameterization by a 3-valent graph $\Gamma$
then one can equip the edges of $\Gamma$ with some extra data
called {\em the phases} that determine $B$.
Let $E$ be a phase of weight $w$
and parallel to a primitive integer vector $v\in\Z^n$.
The vector $v$ determines an equivalence relation $\sim_v$ in the torus $T^n$.
We have $a\sim_v b$ for $a,b\in T^n$ if $a-b$ is proportional to $v$.
Clearly, $T^n/\sim_v$ is an $(n-1)$-dimensional torus.
The {\em phase} of $E$ is a multiset $\Phi=\{phi_1,\dots,\phi_w\}$,
$\phi_j\in T^n/\sim_v$
(recall that $w$ is the weight of $E$).
Alternatively, $\phi_j$ may be viewed as a geodesic circle in $T^n$.
We orient this geodesic by choosing $v$ going away from $A$ along $E$.
A phase determines a collection
of holomorphic cylinders in $\Log^{-1}(E)\subset\tor$.
If some of $\phi_j$ coincide
then some of these cylinders have multiple weight.

Let $A$ be a 3-valent vertex of $\Gamma$ and $E,E',E''$
are the three adjacent edges.
to $A$ with phases $\Phi,\Phi',\Phi''$.
%Suppose that $E,E',E''$ are parallel to the primitive integer vectors
%$v,v',v''\in\Z^n$.
%Let $\sim_A$ be the equivalence relation generated by $\sim_v,\sim{v'}$
%and $\sim_{v''}$. Clearly, $T^n/\sim_A$ is an $(n-2)$-dimensional torus.
%%We have linear maps $T^n/\sim_{v}\to T^n$
The phases are called {\em compatible} at $A$ if
the geodesics of
$\Phi\cup\Phi'\cup\Phi''$ can be divided into subcollections
$\Psi$ such that for every $\Psi=\{\psi_1,\dots,\psi_k\}$
there exists a subtorus $T^2\subset T^n$ containing all geodesics
$\psi_j$ and these (oriented) geodesics bound a region of zero area
in this $T^2$.

\begin{defn}\label{sctc}
{\em A simple complex tropical curve} is a simple tropical curve
$h:\Gamma\to\R^n$ (see Definition \ref{simplecurve}) whose edges
are equipped
with admissible phases such that
%\begin{itemize}
%\item
for every edge $E\subset\Gamma$ the phase
$\Phi=\{\phi_1,\dots,\phi_w\}$ consists of the same geodesic
$\phi_1=\dots=\phi_w$.
%\end{itemize}
\end{defn}

Note that a simple complex tropical curve defines a complex tropical curve
$B\subset\tor$ of the same genus as $h:\Gamma\to\R^n$.
If the phase of a bounded edge of $\Gamma$ consists of distinct geodesics
then the genus of $B$ is strictly greater than that of $C$.

In a similar way one can define {\em real tropical curves} by requiring
all curves $V_t$ in Definition \ref{defB} to be real.
Our next purpose is to define simple real tropical curves.
Let $h:\Gamma\to\R^n$ be a simple tropical curve.
Consider an edge $E\subset\Gamma$ of weight $w$ parallel to a primitive
vector $v\in\Z^n$.
The scalar multiple $wv$ defines an equivalence relation $\sim_{wv}$
in $\Z_2^n$. We have $a\sim_{wv} b$ if $a-b\in\Z_2^n$ is a multiple
of $wv\mod 2$.
The equivalence is trivial if $w$ is even.
Otherwise $\Z_2^n/\sim_{wv}\approx\Z_2^{n-1}$.

The {\em sign} of $E$ is an element $\Z_2^n/\sim_{wv}$.
The choice of signs has to be compatible at the vertices of $\Gamma$.
Let $A$ be a vertex of $\Gamma$ and $E_1,E_2,E_3$ be the adjacent edges
of weight $w_1,w_2,w_3$ parallel to the primitive vectors $v_j\in\Z^n$.
Let $\sigma_j$ be the sign of $E_j$.
%The sign $\sigma_j$ of $E_j$ determines two elements
%$\sigma'_j$ and $\sigma"_j=\sigma'_j+w_jv_j$
%which are distinct if $w$ is odd and equal if $w$ is even.
%Edges $E_j$ and $E_k$, $j,k=1,2,3$ are compatible at $A$ if the equivalence
%classes $\sigma_j$ and $\sigma_k$ have a common element.
We say that the sign choice is {\em compatible} at $A$ if
every element in the equivalence class $\sigma_j$, $j=1,2,3$,
is contained in another equivalence class $\sigma_k$, $k=1,2,3$, $k\neq j$.
\begin{defn}
{\em A simple real tropical curve} is a tropical curve
$h:\Gamma\to\R^n$ whose edges are equipped with signs compatible at every
vertex of $\Gamma$.
\end{defn}

If all edges of $\Gamma$ have weight 1 then this definition agrees
with {\em combinatorial patchworking}, see \cite{IV}.
%The technique described in \ref{secdim2} can be adjusted
%to count simple real tropical curves and
Simple real tropical curves can be used in real enumerative
problems (see \cite{M-en} and \cite{IKS}
for details in the case of $\R^2$).
%In particular,
%it can be used to compute the so-called Welschinger

\begin{figure}[h]
\centerline{\psfig{figure=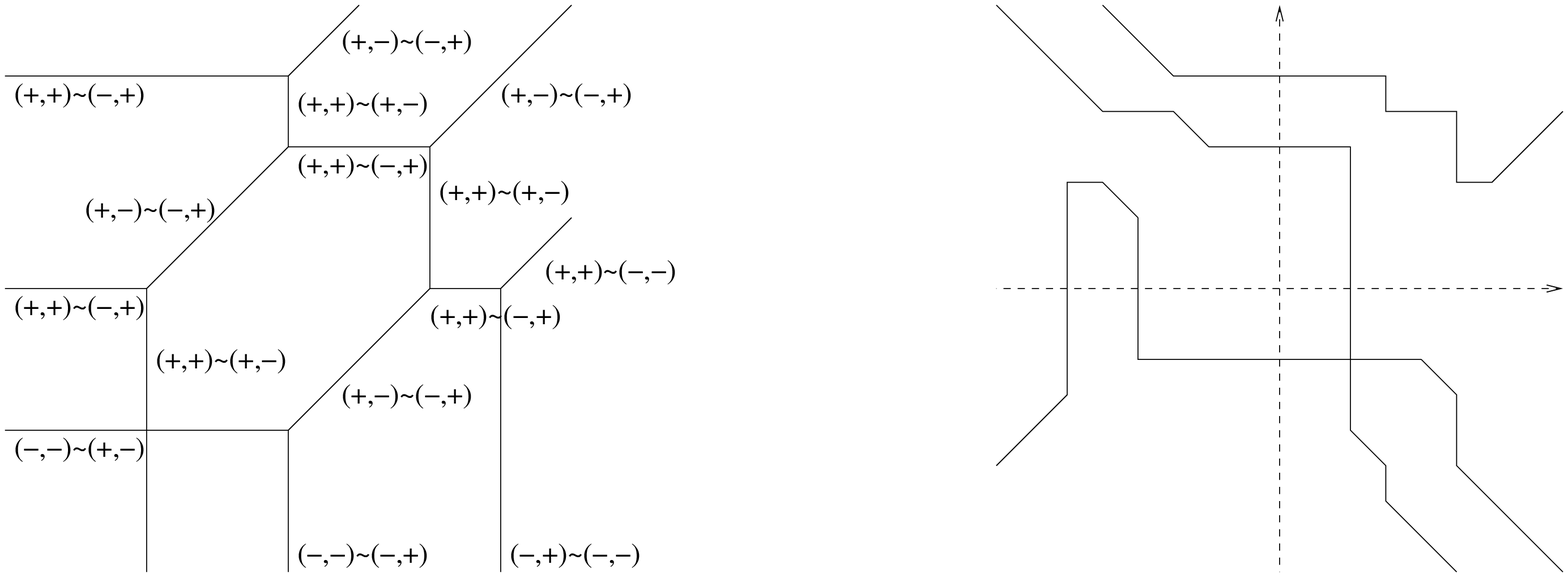,height=1.6in,width=4.4in}}
\caption{\label{realtrop} A real tropical projective cubic curve}
\end{figure}

Figure \ref{realtrop} sketches a tropical curve equipped with
admissible signs and the corresponding real tropical curve.

\end{document}